\overfullrule=0mm
\baselineskip=0,2cm
\lineskiplimit=1pt
\lineskip=1mm
\def \l {{\lambda}}
\def\noi{{\noindent}}
\def\R{{\rm I\kern-2pt R}}
\def\N{{\rm I\kern-2pt N}}
\def\Z{{\rm Z\!\!Z}}

\def\C{{\rm I\kern-5pt C}}
\def\G{{\rm I\kern-5,5pt G}}

\def \e {{\epsilon}}
\def\sp{{\vskip 0,2cm} \noindent }
\def\mp{{\vskip 0,3cm}  \noindent }
\def\bp{{\vskip 0,4cm} \noindent }
\def \noi {{\noindent}}
\def \Ci {{C^\infty}}
 \hsize 16 truecm
\vsize 24 truecm 
\pageno=1
  \centerline{\bf WEIGHTED TRACES ON ALGEBRAS OF
PSEUDO-DIFFERENTIAL OPERATORS}
\bigskip \noindent 
\centerline{ \bf AND GEOMETRY ON LOOP GROUPS}
  \vskip 0,5cm
 \centerline{\it January 2000 }\vskip 0,5cm
 \centerline{\it   A.CARDONA*, C.DUCOURTIOUX,
J.P.MAGNOT, S.PAYCHA } \vskip 0,3cm
 \centerline{Laboratoire de Math\'ematiques Appliqu\'ees}
\centerline{Universit\'e Blaise Pascal (Clermont II)}
\centerline{Complexe Universitaire des C\'ezeaux}
\centerline{63177 Aubi\`ere Cedex}
\vskip 0,2cm
\centerline{(*) and Universidad de Los Andes, Bogota, Colombia}
 \vskip 0,2cm
\centerline{cardona@ucfma.univ-bpclermont.fr}
  \centerline{ c.ducour@ucfma.univ-bpclermont.fr}
 \centerline{
magnot@ucfma.univ-bpclermont.fr}\centerline{
paycha@ucfma.univ-bpclermont.fr}
 \vskip 1cm \centerline{ \bf Acknowledgements} \sp The last author  would like to thank 
 J\"urgen Jost  for 
inviting her for a three months stay at the Max Planck Institut in den  Naturwissenschaften 
during the fall of 1998
 when a big  part  of this 
paper was written. It was completed during a stay 
at the Mathematics Institute of the University of
Bonn supported by the von Humboldt Stiftung. 
  \bp \centerline{\bf Abstract }\vskip 0,5cm
 Using  {\it weighted  traces} which are
linear functionals of the type 
$$A\to  tr^Q(A):=\left(tr(A
Q^{-z})-z^{-1} tr(A
Q^{-z})\right)_{z=0}$$ defined on the whole 
algebra   of (classical) pseudo-differential
operators (P.D.O.s) and where $Q$ is some
positive invertible elliptic operator, we
investigate the geometry of loop groups in the
light of the cohomology of pseudo-differential
operators. We set up a geometric framework  to
study a class of infinite dimensional manifolds in
which we recover some results on the geometry of
loop groups, using again weighted traces.
Along the way, we investigate properties of
extensions of the  Radul and Schwinger cocycles
defined with the help of weighted traces.
\sp   \bp \centerline{\bf R\'esum\'e
}\vskip 0,5cm A l'aide de {\it traces 
pond\'er\'ees}
qui sont des fonctionnelles lin\'eaires du type:  
$$A\to  tr^Q(A):=\left(tr(A
Q^{-z})-z^{-1} tr(A
Q^{-z})\right)_{z=0}$$ d\'efinies sur toute
l'alg\`ebre des op\'erateurs
pseudo-diff\'erentiels classiques, $Q$ \'etant un
op\'erateur elliptique inversible, on \'etudie la
g\'eom\'etrie de l'espace des lacets \`a la
lumi\`ere de la cohomologie des op\'erateurs
pseudo-diff\'erentiels. On met en place un cadre
g\'eom\'etrique  afin d'\'etudier
une classe de vari\'et\'es de dimension infinie,
cadre dans lequel on retrouve, toujours \`a l'aide
des traces pond\'er\'ees, des r\'esultats
concernant la g\'eom\'etrie des lacets. Ces
traces pond\'er\'ees nous permettent aussi
d'\'etendre la notion de cocycle de Radul et de
Schwinger et d'en \'etudier certaines
propri\'et\'es. 
\vfill \eject \noindent
\centerline{\bf Introduction}\bp \bp The paper is
built up from two parts: the first one presents
algebraic tools which are used in the second part
to  extend some geometric concepts to the infinite
dimensional context. The approach to the geometry
of loop groups we present in the second part uses  
weighted traces and it  is new to our knowledge.
Weighted traces offer a useful tool to prove 
both algebraic and geometric results, some of
which had been proved elsewhere by other
methods.     \sp Let us describe the contents of
the first part of the paper (sections 1-4).
 The Lie algebra of interest in our
context is the infinite dimensional
 Lie algebra of
 pseudo-differential operators (P.D.Os) acting on sections of some
finite rank vector bundle $E$ based on a closed manifold $M$ (its dimension
does not yet play a role at this stage).
The Lie bracket is given by the operator bracket.
It is well known that when the manifold is connected and of dimension
strictly larger than $1$,  the only  trace on this algebra, i.e  the only
$\C$ valued linear functional which
satisfies the  {\it tracial   property }  namely $tr[A,B]=0$,
 is the Wodzicki residue trace [W], see also [K] for a review.
\sp Here we shall consider linear functionals   on this algebra which arise
as  {\it zeta function regularized  traces }. They involve a
{\it weight} given by a positive
self-adjoint elliptic operator and in general
depend on the choice of the weight hence the terminology
{\it  weighted  traces } which we shall use here.
One can find such "traces" (sometimes  implicitely)  in    the
literature on determinant bundles [BF] (in particular in the connection and the
curvature) and more generally  when investigating geometry in infinite dimensions (in particular
for the notion of minimality of submanifolds [MRT],[AP]).  Here we will be using them   
to define the Ricci curvature on current groups and first Chern form on loop groups. \sp
 These weighted traces
 extend the so  called canonical trace
of Kontsevich and Vishik  [KV1] defined on  the subalgebra
of P.D.Os with integer order that lie in the odd class (defined below) provided
the dimension
of the underlying manifold is odd.
 On this subalgebra (which contains ordinary differential operators)
the weighted traces
 actually obey the tracial property
$tr[A,B]=0$ which is not the case on the whole algebra of P.D.Os.
\mp Here,
 rather than searching for subalgebras on which
the "weighted traces" actually are traces, we focus
on the obstruction that prevents them from being traces
 on a bigger subalgebra of P.D.Os.
It is measured by the coboundary of the weighted trace and yields
a generalization of  the well known Radul cocycle (see [KK], [R],
  [M2]). It is
of  infinite dimensional essence since it can be
 expressed in terms of a Wodzicki
residue which is a purely infinite  dimensional trace.
   \sp With this idea in   mind of studying
infinite dimensional obstructions rather than  restricting
ourselves to subalgebras on which they vanish,
 using weighted traces
 we build up bilinear functionals on {\it the
whole} algebra of classical P.D.Os. These 
restrict to the   twisted Radul [M2] and Schwinger
[S], [M2], [CFNW] cocycles on the (rather small)
subalgebra of P.D.Os that lie in $g_{res} $, the
algebra of bounded operators
$A$ such that 
$[A,
\e ]$ is 
 of order no larger than
$- { dim M \over 2}$ where $\e$ is the "sign" of
some self-adjoint elliptic operator. \sp Finally,
we investigate the link between our Schwinger
functional and a generalization of another cocycle
arising in the context of a central extension of
the group
$G_{res}  $ of invertible operators in $g_{res}$.
These relations boil down to well-known [PS]
identifications of these cocycles on the
restriction to $g_{res} $ for
  the obstruction to such identifications  which is given in terms of Radul
 cocycles vanishing on this subalgebra.
 \mp Let us now turn to a more geometric
point of view which leads us to the second part
of the paper (sections 5-9). We now consider families of pseudo-differential operators,   extending
the notion of weighted trace to such families. On top of the algebraic obstructions mentioned
above, due to the dependence on the weight, there   are obstructions of
 geometric type which arise when trying to generalize properties
of
 classical geometric concepts to  infinite dimensions. 
\sp Current groups $Map(M, G)$ where $M$ 
is a manifold and $G$ is a Lie group,
 offer a first tractable example of infinite dimensional manifold. As spaces of classical paths 
in Wess-Zumino Witten models, they play an important role in quantum field theories. 
They have been the topic of many an investigation 
and particularly  from a geometric point of 
view (e.g. [DL], [F1,2], [P], [PS], [SP], [Wu]). 
In a pioneering article [F1], Freed   suggested a
way to geeneralize to these groups some  of the
methods available to study finite dimensional Lie 
groups. In particular, he defined a Ricci
curvature on current groups and a first Chern 
form on  the
$H^{1\over 2}$ based  loop group which is
K\"ahler. Other methods [SW] have since then been
suggested, leading  to the same expressions.  In
[F1] (see also [F2], [SW])  the author  uses a
"conditioned" trace which involves taking  a "two
step trace", namely first the trace on the Lie
algebra of the finite dimensional Lie group $G$
and  then, when the operator   obtained this way 
is trace-class, taking its trace. He shows that 
the  (conditioned) 
 first Chern form on the based loop group is proportional to the symplectic form, from which follows on one hand that 
   it is closed and hence defines the first Chern class, and on the other hand that   it 
is   K\"ahler-Einstein.   
  \mp Here we suggest a more general approach to
defining   Ricci curvature and first Chern forms on
 a class 
of manifolds which includes the current  groups mentioned above.    Carried out (in a left invariant
way ) to current  Lie groups equipped with a weight
  given by a left invariant field of elliptic 
operators, the notion of weighted trace  
enables us to define  the {\it weighted   Ricci
curvature} as a weighted pseudo-trace of $R : Z\to
\Omega (Z,\cdot)\cdot$ 
 where $\Omega $ 
is the curvature tensor and the {\it weighted first Chern form} 
as a weighted pseudo-trace of   the curvature,
 provided the operators involved are classical
P.D.Os.
 \sp We express the weighted first Chern form in
terms of a pull-back by the adjoint
representation of  a cocycle on $g_{res}$, thus
relating the closed two form given by the Chern
form with a closed two cochain.  
\sp A well known result by Kuiper [Ku] shows that
a Hilbert manifold is parallelisable since
$GL(H)$ is contractible for a Hilbert space $H$
(the model space of the manifold). This might seem
in conflict with the fact that the first Chern
class does not vanish. But in our approach, the
model space $H$ being a space of sections of some
(finite rank) vector bundle
$E$ based on a closed manifold $M$, the structure group
$GL(H)$ reduces to a non contractible group, namely the
group 
 $Ell_0^*(M, E)$ of zero order invertible elliptic
operators acting on these sections. \mp 
This opens a road to many questions, such as  
finding a criteria  for the weighted first
Chern  form and higher order forms to be
closed, investigating the holonomy group
which we expect to be non trivial because of the 
reasons mentioned above. This article
yields a geometric setting in which   such
questions make sense for a class of Hilbert 
manifolds beyong the example of loop groups.
It confronts this approach with other
approaches in the specific case of loop
groups. \vfill
\eject
\noindent
   \vfill \eject \noindent\bp
\bp  {\bf  1.  Weighted traces  on 
the algebra of
(classical) P.D.Os}\mp {\it In [KV 1], the authors
introduced a new trace  type functional
$TR$ on operators in  $PDO (M, E)$  with order in
$\alpha_0+
\Z$  where
$\alpha_0$ is some {\it non integer} complex
number called the {\it canonical trace}. Avoiding
  integer orders  has to do
with the fact that a positive homogeneous distribution on $\R^n/\{0\}$ of
 non integer  order  can be extended in a unique way to a positive
homogeneous distribution  on $\R^n$ [H].
 In this section we first briefly recall the construction of this trace $TR$ and
some of its properties.  We use this trace as a 
tool to build linear functionals of
pseudo-differential operators  which we call
"weighted traces" and to study some of 
their properties. Such functionals were also
considered in [MN] to prove a pseudo-differential
generalization of the Atiyah-Patodi-Singer index
theorem.  For the sake of
self-containedness, we recall the proof of this
fact.  We refer the reader to Appendix A for
notations and basic facts concerning
pseudo-differential operators. }
\mp
$\bullet$ {\it The Kontsevich and Vishik canonical
 trace}  \mp  Let us first describe the general
lines of the construction in a heuristic way. For
a (classical) pseudo-differential operator locally
given by:
$$Au(x):=\int_{\R^n}   a(x, \xi) \hat u(\xi)d 
\xi $$
where $\sigma_A(x, \xi)$ is the 
(locally defined) total symbol (see Appendix A), we
would like to define $$"TR(A):= \int_M \int_{\R^n} tr_x \sigma_A(x, \xi) d \xi
dvol(x)" $$ which in general does not make sense
since
$a(x,\xi)$ typically has components of degree $\geq- dim M$. It
however does make sense for an operator of order $<-dim M$ and
yields the ordinary trace. We should therefore find a way of only
picking up the finite term in such an expression.
\mp  Let us consider a
 classical P.D.O $A$ of order $\alpha$ with
  symbol   $ \sigma =
\sum_{j=0}^N \Psi  \sigma_{\alpha-j}+ \sigma_{(N)}$ ( see [Sh] for a
discussion about such assumptions on the symbol) where
$\sigma_{\alpha-j}\in S^{\alpha-j}(\R^n)$, 
$\sigma_{(N)}\in S^{\alpha-N-1}(\R^n)$
and where $\Psi$ is a smooth function on $\R^n$
which is zero in $B(0, {1\over 4})$ and equal to
$1$ on $\R^n-B(0, {1\over 2})$. Since $
\sigma_{(N)}(\xi)= O (\vert \xi \vert^{\alpha-N})$
we have $ \int_{B(0,R)} \sigma_{(N)}(\xi) d\xi =
\int_{\R^n} \sigma_{(N)}(\xi) d\xi+
O(R^{\alpha+n-N})$. On the other hand, splitting
the integrals $\int_{B(0, R)}\Psi
\sigma_{\alpha-j}(\xi) d\xi =
\int_{B(0,1)}\Psi \sigma_{\alpha-j}(x,\xi) d\xi + \int_{B(0, R)/B(0,
1)}\Psi \sigma_{\alpha-j}(x,\xi) d\xi $, we use  the fact that
$\sigma_{\alpha-j}$ is homogeneous of order $\alpha-j$ to express the
last integral. If $\alpha$ is an integer then
 there is an integer  $j_0$
such that $\alpha-j_0+n=0$ and we have for $N>j_0$:
$$\eqalign{ &\sum_{j=0}^N \int_{B(0, R)/B(0,
1)}\Psi \sigma_{\alpha-j}(x,\xi) d\xi 
=\sum_{j=0}^N\int_1^R\int_{\vert \xi\vert=1}
r^{\alpha-j+n-1} \sigma_{\alpha-j}(x,\xi) d\xi
dr\cr &\sim_{R\to\infty}
\sum_{j=0,n+\alpha-j\neq 0 }^N{1\over \alpha+n-j}R^{\alpha-j+n}
\int_{\vert \xi\vert=1} a_{\alpha-j }(x,\xi) d\xi+ log R
\int_{\vert \xi\vert=1} a_{-n}(x,\xi)d\xi + \, \hbox{constant
term}. \cr}$$ Finally this yields the existence of an asymptotic
expansion in $R\to +\infty$ and of a constant (w.r.to $R$)   $C(a(x,
\cdot))$ such that $$\eqalign{ &\int_{B(0, R)} a(x,\xi) d\xi\cr &
\sim_{R\to\infty} \sum_{j=0,n+\alpha-j\neq 0 }^N{1\over
\alpha+n-j}R^{\alpha-j+n} \int_{\vert \xi\vert=1}
\sigma_{\alpha-j}(x,\xi) d\xi+ log R \int_{\vert \xi\vert=1}
a_{-n}(x,\xi)d\xi + C(a(x,\cdot)).\cr} $$Because
of the logarithmic term in $R$, one does not
expect   the finite part $$ f.p. \left(\int_{\R^n}
a(x,\xi) d\xi\right)  = LIM_{R\to\infty} \int_{B(0, R)}
\sigma (x,\xi)d\xi  $$ to be invariant
under a change of variable of $\R^n$. However if
the order $\alpha$ is not an integer, there is no
logarithmic divergence and $f.p. \int_{\R^n}
\sigma(x, \xi) d \xi$ is independent of the local
representation $a(x, \xi)$ of
$A$: \mp{\bf
Proposition [KV1, 2] (see also [Le] in Prop. 5.2):} {\it Provided
$A\in PDO(M, E)$ has non integer order:
 $$TR(A):=
\int_M tr_{E_x}\left(f.p. (\int_{\R^n} \sigma_A(x,\xi) d\xi)\right)dvol(x)
\eqno(1.1) $$ where $\sigma_A$ is the symbol of
$A$, is well-defined and satisfies the tracial
property:
$$TR([A,B])=0\quad \forall A\in PDO(M, E), B\in PDO(M, E),
\hbox{such that} \quad  ord(A)+ord(B)\notin \Z. $$ 
Here $tr_{E_x}$ denotes the trace on the fibre $E_x$ of the bundle $E$ above the point $x$.}
\sp The linear functional $TR$
coincides with the usual trace for P.D.Os of
order with real part strictly smaller than minus
the dimension of the underlying manifold the
operators are acting on.
  $TR$ is in fact a trace functional i.e
   $TR[A,B]=0$ for any $A, B\in PDO(M, E)$ such that $ord A+ord B\in \C-\Z$
(see Proposition 3.2 in [KV2]).\sp
    \mp $\bullet$ {\it A fundamental property of the canonical trace}\sp
   \sp Following Kontsevich and Vishik, we shall call
a local family $A_z \in PDO(M, E)$ with distribution kernels (locally)
 {\it weakly holomorphic}   if the following conditons are satisfied:
 \sp (i) The order $\alpha_z$ of $A_z$ is (locally) holomorphic in $z$,
\sp (ii)The kernel
 $A_z(x, y)$ of $A_z$ is (locally) holomorphic in $z$
 for  $x,y$ in disjoint local charts,
\sp (iii)Given any local chart $U$ on $M$,
the homogeneous components $\sigma_{A_z,\alpha_z -j} (x, {\xi\over {\vert
\xi\vert}})$ of the symbol $\sigma_{A_z}(x, \xi)$ are (locally) holomorphic
functions in
$z$ on the restriction to $U$ of the cotangent sphere  bundle $S^*M $.
\sp (iv) When $x$ and $y$ belong to a common local chart $U$, the difference
between $A_z(x,y)$ and the truncated kernel $\sum_{j=0}^N \int \rho(\vert
\xi\vert)  \sigma_{A_z,\alpha_z -j} (x, \xi) exp(i\langle x-y, \xi\rangle d\xi$ of class
$C^{k(N)} $ for some $k(N)$ increasing with $N$  tends  to a (locally)
holomorphic kernel on $U\times U$ when   $N\to \infty$.
\sp In (iii) the topology on the space of symbols is given by the
 supremium norm of the symbol and its derivatives and in (iv) the convergence is to be understood in the 
sense of weak convergence of distributions [KV1]. 
\bp
The following   property of the canonical trace  plays a
fundamental part in these notes:
\mp \noi {\bf Fundamental  property}(see [KV2] Proposition 3.4  and
[KV1] Th.3.13) \sp
{\it For any (local weak) holomorphic family $A(z)$ of classical P.D.Os
on $M$, $z\in U\subset \C$, $ord A_z =\alpha(z)$ where $\alpha$ is
holomorphic and
$\alpha^\prime $ does not vanish, the function $TR(A_z )$ is meromorphic
with no more
than simple poles at $z=m\in U\cap \Z$.}
\mp $\bullet$ {\it Wodzicki residue } \sp Applying
the fundamental property to the family $A^Q_z := A
Q^{-z}$ where
  $Q \in
Ell_{ord >0}^{*, +}  (M, E)$, $A \in PDO(M, E)$ leads to the notion of {\it
Wodzicki residue}  [W] (see also [K] for a review):
$$ res( A  ):= ord Q\cdot Res_{z=0}TR(A^Q_z )
\eqno(1.2) 
$$
  which is in fact independent of $Q$. This can be carried out 
for  any operator  $Q \in Ell^{+}_{ord >0}(M, E)$ which might not be injective,
 replacing it 
 by the operator $Q+P_Q$ in the above formulas so that $A^Q_z= A (Q+P_Q)^{-z}$
   where $P_Q$ denotes the orthogonal 
projection of $Q$ onto its kernel which is finite dimensional since $Q$ 
is elliptic and the manifold closed. The projection is orthogonal for the 
inner product on the space $\Ci(M, E)$ of smooth sections of $E$ 
induced by the hermitian structure $\langle \cdot, \cdot\rangle_x$ on the fibre over $x\in M$ 
 and the Riemannian volume measure $\mu$ 
 on $M$:
$$\langle \sigma, \rho\rangle:= \int_M d \mu(x) \langle \sigma(x), 
\rho(x)\rangle_x\quad \forall \sigma, \rho \in \Ci(M,  E)
$$ 
 \sp 
The Wodzicki residue   defines a trace
functional on the algebra  $PDO (M, E)$ of classical P.D.Os  and
 vanishes for  any classical P.D.O with  non integer order. Since $res$
also vanishes
 on a smoothing operator, it induces a trace functional
 on the  symbol algebra of  $PDO(M, E)$ [W]. It is in fact the unique
trace functional on $PDO(M, E)$ provided $M$ is
  connected and has dimension $>1$.\sp An important feature of the  Wodzicki
 residue we need
 to keep in mind here is that
it  vanishes on any trace-class operator
as can easily be checked from the above
definition.
\mp Kontsevich and Vishik (see [KV2]
  see (3.16))
furthermore show that, given  a local weak
holomorphic family
$A_z$ of operators of order $\alpha(z)$,
the poles of $TR(A_z)$ at entire points
  are expressed in terms of Wodzicki
residues:
  $$Res_{z=m}TR(A_z )=
-{1\over
\alpha^\prime(\alpha^{-1}(m))}
res\left(\sigma(A_{\alpha^-1(m)}
\right).\eqno(1.3)$$
   \mp
 \mp $\bullet$ {\it The weighted trace }\sp
 The
fundamental property leads yet to another linear
functional on the algebra $PDO(M, E)$ which  is
not a trace but interesting all the same because
it does not vanish on trace class operators for
which it coincides with the
 ordinary trace.
\sp  Given  $Q \in Ell^{*, +}_{ord >0}(M, E)$ of order $ord Q$,   we call  the
{\it  $Q$-weighted
trace } of an operator $A\in PDO(M, E)$:
$$tr^Q(A):= \left[TR(A Q^{-z})-
{1\over ord Q\cdot  z }res( A )\right]_{z=0}.
\eqno(1.4)$$ Here again, this extends to the case
when
$Q$ is not injective setting:
$$tr^Q(A):= \left[TR(A (Q+P_Q)^{-z})-
{1\over ord Q\cdot  z }res( A )\right]_{z=0}
\eqno(1.5)
$$ where as before $P_Q$ is the orthogonal
projection onto the kernel of $Q$.  
\sp {\it Warning! To simplify notations, we shall 
 often assume that $Q$ is injective subintending that when it is not, one should 
 replace 
$Q$ by $Q+P_Q$}.  
\sp 
   As we shall soon see, although
refered to as weighted traces here, these functionals do not satisfy the
tracial property $tr^Q[A, B]=0$
and hence do not deserve the name "trace". We shall all the same keep to
 this abusive terminology which turns out to be very convenient. 
\sp $\bullet$ {\it Dependence on the weight}\sp
 Unlike the Wodzicki residue, it  generally  depends on the choice of $Q$.
The 
dependence is  intrinsically infinite dimensional 
since it is
  measured in terms of a Wodzicki residue. Indeed, let  $Q_1, Q_2$
be two
operators in $Ell^{*, +}_{ord >0}(M, E)$ with  same order $q$,  applying the
 fundamental property and $(*)$ to
the (locally around zero) holomorphic family $A({Q_1^{-z}-Q_2^{-z}\over z})$
 of order $\alpha(z):= ord A-zq$, we find (see Prop.2.2 in [KV2]):
$$ \eqalign{ tr^{Q_1}(A)-tr^{Q_2}(A)&=
\lim_{z\to 0} TR(A({Q_1^{-z}-Q_2^{-z}\over z}))\cr
&=  -q^{-1} \cdot res ( A (log  Q_1-log
Q_2))\cr}\eqno(1.6)$$ which is well defined since
$log  Q_1-log Q_2$ lies in $PDO(M, E).$
  \sp However
 when  the underlying manifold is odd dimensional and
 for any odd-class  classical P.D.O $A$ with integer order, the
TR-generalized zeta  function  $z \to TR(A Q^{-z})$ is regular at
$0$ and 
$ tr^Q (A) $ obtained as the limit when $z\to 0$ of these expressions
 is    independent of the
 choice of
$Q$ (see [KV1] Prop. 4.1  where it is   denoted by $Tr_{(-1)}$).
This limit can be seen as an extension of the canonical trace
to operators with integer orders.
\mp 
Although   weighted traces are not tracial, they
have a  useful covariance property:
 \sp {\bf Lemma 1}\sp 
{\it Let $C\in PDO(M, E)$ be injective and bounded. With the same notations as above, we have: 
$$tr^{C^{-1}QC}(A)= tr^Q(CAC^{-1}).\eqno(1.7) $$}
 {\bf Proof:}
 For $z\in \C$ with real part large enough, we have:
 $$\eqalign{ TR(A (C^{-1} Q C)^{-z})&=  TR(A C^{-1} Q^{-z} C) \cr 
&=  TR(   C A C^{-1} Q^{-{z }}) \cr}$$ 
where we have used that $TR$ is  tracial.
 Taking the renormalized limit then  yields the
result.\hfill $\bullet$ 
\bp {\bf 2.  The Radul cocycle as a coboundary of 
the weighted trace}\mp
 {\it In this section we investigate the coboundary of
 the weighted traces introduced
 above. We show how the  cocycles  obtained in this way  relate  to the
  Radul cocycle which arises  in
 geometric quantization [M2]. In the context of
regularized determinants, it
is related to the multiplicative anomaly [D].
This cocycle was already investigated in [MN]
where  the authors    express 
 the coboundary of weighted traces in terms of a
Wodzicki residue (see [MN] Lemma 13). } \mp
We shall need some  definitions of Lie algebra cohomology (see e.g. [M1] ).
\mp $\bullet$ {\it Lie algebra cohomology:} Let $L$ be a Lie
 algebra and $V$ an $L$-module (the action of $L$ on $V$ is denoted by a "dot").
 A {\it cochain of degree n} ( or {\it n-cochain}) with values in $V$
is an antisymmetric multilinear map
$c: L\times L\times \cdots \times L $ (n times) $\to V$. Let $C^n(L, V)$ denote
 the space of all $n$-cochains and let us define the coboundary operator:
$$ \delta: C^n(L, V)  \to C^{n+1} (L, V)$$
$$\eqalign{ (\delta c^n) (x_1, x_2, \cdots, x_{n+1}) =
&\sum_{i<j} (-1)^{i+j} c^n([x_i, x_j], x_1,\cdots, \hat x_i, \cdots, \hat
x_j, \cdots, x_{n+1})\cr
&+ \sum_{i=1}^{n+1} (-1)^{i+1}x_i\cdot c^n(x_1, \cdots, x_{i-1}, x_{i+1},
\cdots, x_{n+1})\cr}$$
Taking $V=\C$ with  the trivial zero action  of $L$ on $\C$, $x\cdot z :=0$
 for any $x\in L $, $z\in \C$  we have :
$$ (\delta c^n) (x_1, x_2, \cdots, x_{n+1}) =
 \sum_{i<j} (-1)^{i+j} c^n([x_i, x_j], x_1,\cdots, \hat x_i, \cdots, \hat x_j,
 \cdots, x_{n+1}).$$
In particular $\delta^2=0$ . An {\it n-cocycle} is a
cochain of degree $n$
 with vanishing coboundary. Let us
   denote by $B^n(L, V)$ the set  of $n$-coboundaries, by
$Z^n(L, V)$ the set of $n$ cocycles and let us call  $H^n(L, V):= Z^n(L,
V)/ B^n(L, V)$
the $n$-th   cohomology space.
\mp
$\bullet$ {\it The weighted  Radul cocycle on 
 $PDO(M, E)$}\sp 
We  now  apply this construction to $L:= PDO(M, E)$ and the $1$-cochain
given by a weighted trace.
 \mp   Let $Q \in Ell^{*,+}_{ord >0}(M, E)$. The $Q$-weighted traces
$tr^Q$ do not satisfy the cyclicity property
thus leading to
 a cocycle  given by its coboundary:
$$c_R^Q(A, B):= \delta tr^Q( A,B )=tr^Q[A,B]\quad 
\forall A, B\in PDO(M, E) \eqno(2.1) $$
which we call the $Q$-weighted Radul cocycle. This terminology is justified
on the  grounds of
the following proposition  which relates $c_R^Q$ to a more familiar
 expression of the classical Radul cocycle.
 \mp {\bf Proposition 1}\sp  {\it (i) For any $A  \in PDO(M, E)$  the
operator $[log Q    ,  A ]   $ lies in
$PDO(M, E)$. \sp (ii) For any $A,B\in PDO(M, E)$ we have:
 $$c_R^Q(A,B)= - {1\over ord Q}res(  [log Q    ,  A ] B ) =
\lim_{z \to 0}TR(Q^{-z} [A,B]).\eqno(2.2)$$   }
\mp {\it Remark}\sp
{\it This relation expresses the algebraic obstruction preventing
a weighted trace from being tracial
in terms of a Wodzicki residue, which is a trace of purely   infinite
dimensional type.}
 \mp   {\it Proof}: \sp 
 (i) We first  check    that
 $[log Q    ,  A ]   $ lies in $PDO(M, E)$ if the order $ord A$ of $A$ is
strictly positive.
$$\eqalign{ [A, log Q ]
&=(A \cdot log Q-log Q\cdot  A)   \cr
&= A\left(log Q- {ord Q\over ord A}log (\vert A\vert  +1)\right)+
\left({ord Q  \over ord A} log (\vert A\vert  +1)-log Q\right) A . \cr} $$
 Since  the difference of two logarithms of PDOs of same order is a P.D.O,
this proves that if the order of $A$ is strictly positive then  $[ log Q ,
A]   $ is  a P.D.O.
Now if $A$ has negative order,  let $k=  {\vert ord A\vert +1\over ord Q}
$. Then setting $P= Q^k$
we have that $P$ is  of order equal to $\vert ord A\vert +1$ and  we have:
$$\eqalign{ P(Alog Q-(log Q )A)&= Q^k A \  log Q-log Q\  Q^k A\cr
&= [Q^k A, log Q],\cr}$$
  which by the previous results applied to $Q^k A$ which has strictly positive
order,  shows it is a P.D.O.  \sp
(ii) For $Re z$ large enough we have:
$$ \eqalign{TR(Q^{-z} [A,B])&=TR([Q^{-z}, A] B) + TR ( [A, Q^{-z } B])\cr
&=TR([Q^{-z}, A] B) +TR([A Q^{-z\over 2}, Q^{-z\over 2} B])\cr &=
TR([Q^{-z}, A] B).\cr}$$
In the same way we have:
$$ \eqalign{TR(Q^{-z} [A,B])&=TR([A, B]Q^{-z})\cr
&=TR(A[B,Q^{- z}]  )+ TR( A Q^{-z} B-BA Q^{-z})   \cr
&=TR(A[B,Q^{- z}]  ) +TR([ A Q^{-z\over 2}, Q^{-z\over 2} B ])\cr &=
TR(A[B,Q^{- z}]  ).\cr}$$
 Appplying the fundamental property and $(*)$ to the family
$z^{-1}[Q^{-z}, A] B$ of order  $\alpha(z)= -  z\cdot ord Q+ord A+ord B$
we find:
$$\eqalign{ \lim_{z\to 0}TR ( [Q^{-z}, A] B )&=Res_{z= 0}
TR ( z^{-1}[Q^{-z}, A] B )\cr
&= {1\over ord Q}res\left(  {d\over dz}_{/z=0}( [Q^{-z}, A] B  )\right)\cr
&=-{1\over ord Q} res(  [log Q  , A] B) .\cr}$$
 In a similar way we have:
$$ \lim_{z\to 0}TR(A[B,Q^{- z}])=-{1\over ord Q} res(   A[B, log Q  ]
 ).$$
\sp  This proposition shows that the
 $Q$-weighted  Radul  cocycle    generalizes the usual Radul cocycle
obtained for  $\sigma(Q)=\vert\xi\vert$  [CFNW] (see formula (24)),[M2],
[R]  (see formula (41)), [KK]  (see formula (1)).
In this particular case, the Radul cocycle 
 was also considered  in [KV2] in
relation to multiplicative anomalies for determinants of elliptic
operators (see also [D]) but it was later 
considered for more general operators $Q$  in [MN] 
and proved (see Lemma 13) to be  a coboundary  in
the Hochschild cohomology of pseudo-differential
operators. As a consequence of the above
proposition,
 the weighted Radul cocycle vanishes on the algebra of odd-class P.D.Os with
integer order whenever the underlying manifold is odd dimensional.
\mp To each $\C$-valued cocycle $c^Q_R$ corresponds a central extension
(see e.g [M1], [Ki], [Rog])
 which we shall denote by
 $ PDO(M,E)^Q:= \{ (A, \l), A\in PDO(M, E)\}$
 with Lie bracket:
 $$ \left[ (A, \l), (B, \mu)\right]^Q:= \left( [A, B], c_R^Q(A,
B)\right).$$ In other words we have the exact sequence of Lie algebras:
 $$0\rightarrow \C \to PDO(M,E)^Q\to PDO(M,E)\to 0$$
 Two such extensions $PDO(M, E)^{Q_1}$ and $PDO(M,E)^{Q_2}$ are
equivalent since for
 $Q_1, Q_2\in Ell^{+*}_{>0} (M, V)$ the cocycles $c_R^{Q_1}$ and
$c_R^{Q_2}$ are cohomologous. Indeed  their difference
    $$c_R^{Q_1}-c_R^{Q_2}=res\left( (log Q_1-log Q_2)[\cdot,
\cdot])\right)$$ is the coboundary of the $1$ cochain
 $A\to res\left( (log Q_1-log
Q_2)\quad \cdot\right)$.
 \sp
 Since the difference of two logarithms is a zero order P.D.O, this latter
$1$ cochain is a particular example of a   family  of cochains parametrized
by the
 algebra $PDO^0(M, E)$ of operators in $PDO(M, E)$ of order $zero$:
 $$res^P:= res(P\quad \cdot) \eqno(2.3)$$
 where $P$ is a zero order P.D.O. The coboundary
$(A, B) \to res(P[A, B])$ does not vanish in
general. The case when $P=\e(D)$ is the sign
of a self-adjoint operator $D$ gives rise to 
Mickelsson's  $res^\prime$  linear functional.
  \bp {\bf 3. A weighted Schwinger functional}
\mp In this section $D$ denotes a   self-adjoint elliptic operator acting
on smooth sections of a vector bundle $E$ with strictly positive order.
An  example of such a bundle is given by   the spinor
 bundle over   an odd dimensional spin
compact manifold $M$ and the operator $D$ by  the Dirac operator  on $M$.
\sp
Adopting notations which are   frequently used in
the context of geometric quantization,
let
$\e(D):=  D+P_D  (\vert D\vert +P_D)^{-1}$ 
denote the
 sign of $D$ which defines 
 a classical P.D.O. of order $0$.
The operator  $Q:= \vert D\vert $ lies in $  Ell^{+ }_{ord>0} (M, E)$.
 \mp $\bullet$ {\it A Schwinger   functional} \sp
We define the Schwinger
 functional:
$$\eqalign{PDO(M, E)&\to \C\cr
(A, B)&\to  c_S^D(A, B):= {1\over 2}tr^{\vert D\vert} (\e(D) [\e(D), A]
[\e(D), B]).\cr}\eqno(3.1)
$$
The terminology "Schwinger functional" is motivated by the fact 
(as we shall see later) that 
 on some subalgebras of $PDO(M, E)$ it coincides with the
usual Schwinger cocycle.
\sp
A straightforward  computation yields:
$$c_S^D(A, B)= tr^{\vert D\vert} ([A, \e(D)] B)=  tr^{\vert D\vert} ( A[
\e(D), B]).$$
\sp  Let us   introduce the  linear functional
$$\eqalign{ PDO(M, E)&\to \C\cr
A &\to tr^D_\e(A):= tr^{\vert D\vert} (\e(D) A)\cr}$$
which we call the {\it signed weighted trace} of $A$.  As before the
terminology "trace" is not appropriate here since this functional does not
in general have vanishing coboundary.
  \mp   {\bf Proposition 3}\sp {\it  Let ${\cal S}$ be a subalgebra of
$PDO(M, E)$.
The following conditions are equivalent (all the cocycles are Lie algebra
 cocycles):
\item{1)} $c_S^D$ is a $2$-cocycle on   ${\cal S}$,
\item{2)} $c_S^D$ is an antisymmetric 
 bilinear form on   ${\cal S}$,
 \item{3)}The bilinear map $$(A,B)\to 
c_{TR}^D(A,B):= tr^{\vert D\vert}
[\e(D) A, B]$$ defines
a $2$-cocycle on ${\cal S}$,
\item{4)}$c_{TR}^D$ is an antisymmetric  bilinear form  on $S$,
  \item{5)}  For any $ A, B\in {\cal S}$ we have:
 $$  res^{\e(D)}[A, [log \vert D\vert, B]]= 0,
\eqno (3.2) 
$$
  \sp   Provided one of these conditions is
fufilled, the following relation holds:
$$  c_{TR}^D  -  c_{S}^D =
    \delta tr^D_\e   $$
 so that the cocycles $c_{TR}^D$ and $c_S^D$ are cohomologous.
 }\mp
{\it Remarks }\sp  1)  $c_{TR}^D(A,B)$
 coincides with $ -2 c^\prime(A,B)$ defined in (2.11) of  [M2].
\sp 2) Part 5) of this proposition shows that the obstruction
to the cocycle property of the various functionals involved arises as a
 residue
 $ res^{\e(D)}[A, [log \vert D\vert, B]]$.
 Here again the obstruction is therefore
 purely infinite dimensional.   
\mp
 The following definitions and the lemma below will be used in the proof.\sp
We shall set:
  $$\tilde c_{TR}(A,B):= tr^{\vert D\vert} [ A\e(D), B]$$
and $$\bar c_{TR}  :={c_{TR}^D+ \tilde c_{TR}^D\over 2}.$$
\mp
{\bf Lemma 2}\sp
{\it
$$c_{TR}^D(B, A)= -\tilde c_{TR}^D(A,B)\quad \forall A, B\in PDO(M, E)$$ so that
 $\bar c_{TR}^D$ is antisymmetric. Moreover
$$\eqalign{ \bar c_{TR}^D(A,B)&= \delta tr_\e^D (A,B) +{1\over 2}\left(
tr^{\vert D\vert}
([\e(D),B] A) - Tr^{\vert D\vert}([A, \e(D)] B )\right)\cr
&= {1\over 2} \delta tr_\e^D (A,B) +{1\over 2}tr^{\vert D\vert }
\left(  A\e(D) B  -    B\e(D) A  \right)\cr}$$
where $ \delta tr_\e^D $ denotes the coboundary of the signed weighted
trace in the Lie algebra
cohomology.}
\sp \noi {\it Proof}: Let us first see how $c_{TR}^{ D }$ transforms when
exchanging $A$ and $B$.
$$\eqalign{  c_{TR}^D(A, B)&= tr^{\vert D\vert}  [\e(D) B, A]\cr
&= tr^{\vert D\vert} (\e(D) BA- A \e(D) B)\cr
&= -tr^{\vert D\vert} (A\e(D) B-BA\e(D))\cr
&=- tr^{\vert D\vert} [A\e(D), B]\cr
&= -\tilde c_{TR}^D(A, B)\cr}$$
 From this it follows that $\bar c^D_{TR}$ is antisymmetric. The second
statement then easily follows from
Lemma 2 and the definition of $\bar c_{TR}^D$:
$$\eqalign{ \bar c_{TR}^D(A, B)&= {1\over 2} \left( c_{TR}^D(A, B)
 -c_{TR}^D (B, A)\right)\cr
&=  {1\over 2} \left(c_\e (A, B)- c_\e(B, A)\right)+ {1\over 2}
\left(tr^{\vert D\vert}([\e(D), B] A- tr^{\vert D\vert}([\e(D), A] B)\right)\cr
&= c_\e^D(A, B)-  {1\over 2} \left(tr^{\vert D\vert}([\e(D), B] A+
tr^{\vert D\vert}([A, \e(D)] B )\right)\cr
&= {1\over 2}c_\e^D(A, B)+  {1\over 2} \left(tr^{\vert D\vert}( A\e(D) B-
B \e(D) A)\right).\cr}$$
    \mp
   {\it Proof of  proposition 3}:\sp
  $2) \Leftrightarrow  4):$  A straightforward
computation yields:
$$\eqalign{c_S^D(A, B)&=-c_S^D(B, A)\cr
&\Leftrightarrow tr^{\vert D\vert}(A[\e(D), B]) =-tr^{\vert D\vert}
(B [\e(D), A])\cr &\Leftrightarrow
tr^{\vert D\vert} [A\e(D), B]= tr^{\vert D\vert}[A, B\e(D)]\cr
 &\Leftrightarrow \tilde c_{TR}^D (A,B)=c_{TR}(A,B)\cr
 &\Leftrightarrow  - c_{TR}^D (B,A)=c_{TR}(A,B)\cr}$$
where we have used the result of Lemma 2.
\sp  $3) \Leftrightarrow  4):$
The implication from left to right is clear since a $2$-cocycle is
antisymmetric. Let us prove the other implication $4)\Rightarrow 3)$  which
amounts to
 showing that  when $c_{TR}$ is antisymmetric  it defines  a cocycle.
 Since $c_{TR}^D$ is antisymmetric, we have $c_{TR}^D=\tilde  c_{TR}^D=\bar
 c_{TR}^D$. By Lemma 2, it is
sufficient to show that $\omega$ defined by  $\omega(A,B)\equiv Tr^{\vert
D\vert} (A\e(D) B- B\e(D) A)$ has vanishing coboundary  $\delta \omega$.  A
direct computation  yields:
$$\delta \omega(A, B,C)= c_{TR}^D( A, [B\e(D), C\e(D)])+ c_{TR}^D( B,
[C\e(D), A\e(D)])+ c_{TR}^D( C, [A\e(D), B\e(D)])$$
 from which follows that:
$$\eqalign{ \delta \omega(A\e(D) , B\e(D) ,C\e(D) )&= c_{TR}^D( A\e(D), [B
, C ])+ c_{TR}^D( B\e(D), [C , A ])+ c_{TR}^D( C\e(D), [A  , B  ])\cr
&=\tilde  c_{TR}^D( A\e(D), [B , C ])+ \tilde c_{TR}^D( B\e(D), [C , A ])+
\tilde c_{TR}^D( C\e(D), [A  , B
])\cr&= tr^{\vert D\vert} [A, [B, C]]+  Tr^{\vert D\vert} [B, [C, A]]
+Tr^{\vert D\vert} [C, [A, B]]\cr
&=0\cr}$$
 Since any P.D.O  $A$ can be written  $A= A_1 \e(D)$ where $A_1\equiv
A\e(D)$ is a P.D.O, the result follows.
\sp
 $4) \Leftrightarrow  5):$  Since this computation is   similar to the one
 used in the proof of Lemma 2  we shall skip some
intermediate steps here.
 Let $d$ denote the order of $D$ and let us set $l:= log \vert D\vert$.
 $$\eqalign{  - d \cdot c_{TR}^D(B,A)&=- d \cdot c_R^{\vert D\vert
}(\e(D)B, A)\cr
&= res( \e(D) [l,  B]  A )  \cr
&= res( \e(D) A[l,  B]    ) -  res( \e(D) [A,[l,  B]]    )\cr
 &= -d\cdot \tilde c_{TR}^{  D  }( B, A) -  res( \e(D) [A,[l,  B]]    )   \cr
&=   d \cdot c_{TR}^D(A,B)  -  res( \e(D) [A,[l,  B]]    ) \cr}$$
hence $$ c_{TR}^D(B,A)+  c_{TR}^D( A,B)=-{1\over d}  \cdot
 res(\sigma(\e(D) [A,[l,  B]]   )) $$  and $c_{TR}$ is antisymmetric if and
only if condition  (3.2) is satisfied.\sp \noi
so that $c_{TR}^D$ is antisymmetric if and only if $ res( \e(D) [A,[l,  B]]
)=0$
for any $A, B$ in the algebra under consideration.
\sp
$1) \Leftrightarrow  2):$ Only the implication from right to left is
non trivial. To prove it, we use Lemma 2 once
again by which we have:
$$\bar c_{TR}^D(A, B)= {1\over 2} \delta tr_\e^D(A, B) + {1\over 2} Tr^{\vert
D\vert} \left(A\e(D) B- B \e(D) A\right). $$
On the other hand, since $c_S^D$ is antisymmetric, we have:
$$c_S^D(A, B)= {1\over 4}  tr^{\vert D\vert} (\e(D) \left[[\e(D), A] ,
[\e(D), B]\right] .$$
 Then   a direct computation yields:
$$ c_{S}^D(A, B)= -{1\over 2} (\delta tr_\e^D)(A, B) +
{1\over 2} tr^{\vert D\vert} \left(A\e(D) B- B \e(D) A\right). $$
Hence
$$\bar c_{TR}^D  -  c_{S}^D =
    \delta tr^D_\e   $$
 \sp From this  identity follows that, provided $c_S^D$ and $c_{TR}^D$
 are cocycles,
 then they 
 are   cohomologous.
\mp  $\bullet$ {\it The algebra $PDO(M, E)_{res}^D$} \sp
Let us introduce  a subalgebra $PDO(M, E)_{res}^D$of $PDO(M, E)$  the
definition of  which
is close to the algebra
$g_{res}$ (they coincide up to the fact that we require the operators to
 be P.D.Os) which  plays a  substancial part
in geometric quantization techniques (see e.g. [PS]):
$$\eqalign{ PDO(M, E)_{res}^D&= \{A\in PDO(M, E),
\quad  ord A \leq 0\quad  ord ([A, \e(D) ])< -{dim M\over 2}\}\cr
&= g_{res}^D\cap PDO(M, E)\cr}\eqno( 3.3)$$
 where as before $dim M$ is the dimension of the underlying manifold and
where $ord A$ denotes the order of the
 operator $A$. Here $g_{res}^D:= \{ A \in {\cal B} (L^2(M, E)), \quad 
[A, \e(D) ] \quad \hbox{ is Hilbert Schmidt}\}$ where $L^2(M, E)$ is the closure of 
the space $\Ci(M, E)$ of smooth sections of $E$ for the $L^2$ inner product 
induced by the hermitian structure on $E$ and the Riemannian 
volume measure on $M$ as above, ${\cal B} (L^2(M, E))$ denoting the algebra of bounded operators on this Hilbert space.
 \mp
An immediate consequence of Proposition 3 is the following Corollary:
\sp {\bf Corollary 1}\sp {\it  Let ${\cal S}$ be a subalgebra of $PDO(M,
E)$. If  ${\cal S}$ is stable under the map
 $A\to [log \vert D\vert, A]$ i.e  if
 $$A\in {\cal S} \Rightarrow [log \vert D\vert ,
A]\in {\cal S}\eqno(3.4)$$  and if moreover
 $$ res (\e(D) [A,B])=0 \quad \forall A,B \in
 {\cal S}\eqno(3.5)$$  hold on
the subalgebra $S$,
 then   so does relation (R) and   $ c_{TR}^D $ and
$c_{S}^D$ define   cohomologous cocycles on ${\cal S}$.  }
 \mp \noi{\bf Corollary 2}\sp  {\it
On the algebra $PDO(M, E)_{res}^D$ the  Schwinger functional coincides with
the usual
Schwinger cocycle and $c_{TR}^D$ with the usual twisted Radul cocycle and
we have:
$$c_S^D \equiv    c_{TR}^D. $$
}\sp\noi  {\it Proof}:   We   need to  check
 that $PDO(M, E)_{res}^D$ fulfills assumptions
(3.4) and (3.5) of Corollary 1. As before we set
$l:= log \vert D\vert$. Since   $[\e(D), [l, A]]=
[l, [\e(D), A]]$ , the order of $[\e(D), [l, A]]$
is the same as that of
 $[\e(D), A]$ so that  if $A$ lies in $PDO(M, E)_{res}^D$, so does $ [l,
A]$.  Hence assumption
(3.4) is satisfied on $PDO(M, E)_{res}^D$.
 An easy computation yields:
$$ res( \e(D)[A,B] )=  -{1\over 2} res(  \e(D)[ [\e(D), A], [\e(D), B]]
)$$ which vanishes on
$PDO(M, E)_{res}^D$ since it is the trace of an operator of degree
strictly smaller than $-dim M$, the operator    $\e(D) $ being of order $0$
and the operators   $[\e(D), A]$,  $[\e(D), B]$ being both of order
strictly smaller than $-{dim M\over 2}$.
Hence assumption (ii) is satisfied on $PDO(M, E)_{res}^D$.
   Corollary 1 implies that
 both $c_S^D$ and $c_{TR}^D$ are cocycles
 on $PDO(M, E)_{res}^D$ and that they are
cohomologous.
\sp On $PDO(M, E)_{res}^D$ the cocycle $c_S^D$ reads:
$$c_S^D (A, B)={1\over 2} tr( \e(D) [\e(D), A] [\e(D), B]) $$ where now $tr$ is an ordinary trace
and it therefore coincides  with the usual Schwinger cocycle see e.g
 (according to the author, the
definition might change  by  a constant factor)
 [CFNW], [PS], [M2], [S].
   \sp The expression of $c_{TR}^D$
in terms of a residue  obtained in section 1
(Proposition 1) yields
$$c_{TR}^D(A, B)= -{1\over ord D} res( [log \vert D\vert, \e(D) A] B) $$
 thus relating
$c_{TR}^D$  with Mickelssons [M2] twisted Radul cocycle
(they coincide up to a factor $-2$). \bp
 \bp\noi {\bf 4. Cocycles and group
representations}\mp
 Given $D \in Ell^{s.a}_{ord >0}(M, E)$, we have a natural polarization of the space
$H\equiv L^2(M, E)$ given by:    $$H=H_+(D)\oplus H_-(D)$$   where
$H_+(D)= \pi_+(D)(H)$, $H_-(D)= \pi_-(D)(H)$, $\pi_+(D)={1+\e(D)\over 2}$,
$\pi_-(D)={ -\e(D)+Id\over 2}$. Notice that with this choice for $\e(D)$,
$\pi_+(D)$ is $1$ on $Ker D$ and
$\pi_-(D)$ vanishes on $Ker D$. Of course we could have chosen the other
convention  namely  $\pi_-(D)$ to be  $-1$ on $Ker D$ and
$\pi_+(D)$ to vanish  on $Ker D$.  In this polarization, let us write an
operator $A\in PDO(M, E)$
as a matrix:
$$A\equiv \left[\matrix{A_{++}& A_{+-}\cr
                         A_{-+}& A_{--} \cr}\right].$$
Since $ \e(D)= \left[\matrix{1& 0\cr
                         0&-1 \cr}\right]$, we have  $$[\e(D), A]=   2
\left[\matrix{ 0& A_{+-}\cr
                         -A_{-+}&0 \cr}\right].$$
The operator  $\e(D)$ being a P.D.O, the operator $\left[\matrix{ 0& A_{+-}\cr
                         -A_{-+}&0 \cr}\right]$ is also  P.D.O. \sp \noi
For   $D \in Ell_{ord >0}^{s.a}(M, E)$ we  introduce the following bilinear functional:
 $$\eqalign{ \l^D: PDO(M, E) &\to \C\cr
(A, B) &\to  tr^{\vert D\vert}([A_{++},
 B_{++}]-[A,B]_{++}).\cr}\eqno(4.1)$$
  On $PDO(M, E)_{res}^D$ it coincides with  the Lie algebra cocycle
corresponding
to a central extension of the group:
$$G_{res}^D\equiv \{ A\in GL(H), A_{+-} \quad \hbox{ and }
\quad A_{-+} \quad \hbox{ are Hilbert-Schmidt}\}
\eqno(4.2)  $$
as described in [PS] (6.6.5).
\sp The Schwinger functional $c_S^D$ not being antisymmetric in general
(since according to proposition 3
it then becomes
 a cocycle), it is natural to introduce a {\it
mean Schwinger functional}  in a similar way to what we we did for the
twisted Radul cocycle:
$$\bar c_S^D(A,B):= {c_S^D(A, B)- c_S^D(B, A)\over 2}={1\over 4}
tr^{\vert D\vert}
\left[ [\e(D), A], [\e(D), B]\right].$$
  \mp \noi {\bf Lemma 3}\sp {\it Let $D \in
Ell_{ord >0}^{s.a}(M, E)$,  let $A  \in PDO(M, E)$,
$B
\in PDO(M, E)$. Then
$$ \eqalign{  2 \l^D(A,B)
&=2tr^{\vert D\vert} ( B_{+-}A_{ -+} - A_{+-} B_{-+})  \cr
&=   \bar c_S^D(A, B) -c_R^{\vert D\vert} (A_{+-},B_{-+})-  c_R^{\vert
D\vert} (A_{ -+},B_{ +-}) \cr}$$
Furthermore $\l^D$ is  a cocycle on $PDO(M, E)_{res}^D$ and  the following
relation holds on $PDO(M, E)_{res}^D$:
 $$\eqalign{   {1\over 2}c_S^D(A,B)  &=   \l^D(A,B)\cr
&=TR ( A_{-+}B_{+-} - A_{+-} B_{-+}) \cr
&= {1\over 4} TR\left( \e(D)  [\e(D), A]  [\e(D), B]  \right)\cr}$$
 confirming formula (6.6.6) of [PS]. }
\sp \noi {\it Remark }:  Although the individual components $A_{++},A_{+-},
A_{-+}, A_{--}$ of the two by two matrix of $A$ in the polar decomposition
might not be P.D.Os,
$c_R^{\vert D\vert} (A_{+-},B_{-+})= tr^{\vert D\vert} [A_{+-},B_{-+}]$ and
 $c_R^{ \vert D\vert} (A_{ -+},B_{ +-}) =
tr^{\vert D\vert} [A_{-+},B_{ +-}]$ are well defined since these are
weighted traces of   operator brackets  which are   P.D.Os as projections
of operator
brackets of   P.D.Os of the type  $[[\e(D), A], [\e(D), B]]$.
\sp \noi {\it Proof}: Let us start with the general case. The first
equality follows from a direct  computation  and  the definition of the
weighted Radul cocycle. The second equality
follows from the  fact that
$$\e(D)[ [\e(D), A], [\e(D), B]]= -4\left[ \matrix { A_{+-} B_{-+} - B_{+-}
A_{-+} & 0 \cr
0 & B_{ -+} A_{  +-}- A_{ -+} B_{ +-}  \cr}\right]$$ since this yields:
$$\bar c_S^D
(A,B)=2tr^{\vert D\vert} (  A_{+-} B_{-+}-B_{+-}A_{ -+}  )+  c_R^{\vert
D\vert} (B_{+-},A_{-+})+ c_R^{\vert D\vert} (B_{ -+},A_{ +-}).$$
When considering the case of $PDO(M, E)_{res}^D$, the Radul cocycles $
c_R^{\vert D\vert} (B_{+-},A_{-+})$ and $c_R^{\vert D\vert} (B_{ -+},A_{
+-})$ vanish since
$B_{+-},A_{-+}, B_{ -+},A_{ +-}$ are Hilbert-Schmidt and the weighted
traces $tr^{\vert D\vert}$ become ordinary traces $TR$. This completes the
proof.
\bp \noi Noting  that  $[\e(D), A]=-[\e(D), A]^*$ is equivalent to $  A_{
+-}  =- A_{-+}^* $ for an operator $C:H_+(D)\to H_-(D)$ we shall set:
$$j(C) \equiv \left[\matrix{0&- C^*\cr
                         C& 0 \cr}\right].$$
  \mp \noi {\bf Lemma 4 }\sp {\it Let $D \in
Ell_{ord >0}^{s.a}(M, E)$,  let
$A:H_+(D)\to H_-(D)$, $B:H_+(D)\to H_-(D)$ such that $j(A)  \in PDO(M, E)$,
and $j( B)
\in PDO(M, E)$. Then
$$  \eqalign{  \bar c_S^D(j(A),j(B))  &=\delta (tr_\e^D)(j(A), j(B))\cr
&=  2tr^{\vert D\vert} (B^*A-A^*B)+ c_R^{\vert D\vert}(A, B^*)+
  c_R^{\vert D\vert}(  A^*,B)\cr}$$}
\sp \noi {\it Proof}: Since $[\e(D), j(A)]=-2j(A)$ we have:
$$\eqalign{ \bar c_S^D(j(A), j(B))&= {1\over 4} tr_\e^D([[\e(D), j(A)],
[\e(D), j(B)]]\cr
 &= tr_\e^D([j(A), j(B)])\cr
&= c_\e^D(j(A), j(B))\cr
&= tr_\e^D \left[\matrix{ -A^*B+ B^*A &0 \cr
0& -AB^*+ BA^*\cr}\right]\cr
&= 2 tr^D (B^*A-A^*B) + c_R^{\vert D\vert} (A,B^*)- c_R^{\vert D\vert}
(B,A^*)\cr
&= 2 tr^D (B^*A-A^*B) + c_R^{\vert D\vert} (A,B^*)+c_R^{\vert D\vert} (
A^*, B)\cr}$$
\mp \noi On the grounds of this proposition we set for two operators $A,B$
such that $A^*B\in PDO (M, E)$:
  $$\omega^{D} (A,B)\equiv  -i  tr^{\vert D\vert}(A^*B-B^*A).$$
$\omega^D$ relates to the mean Schwinger functional $\bar c_S$ as
follows:\mp \noi
 {\bf Corollary 3}\sp  {\it Let   $A: H_+(D)\to
H_-(D)$ and  $B: H_+(D)\to H_-(D)$ be  operators
such that $j(A)$ and $j(B)$ lie  in  $PDO  (M,
E)$. Then
 $$ \bar c_S^D (j(A),j(B))=2i\omega^D( A ,  B ) + c_R^{\vert D\vert} (A,
B^*)+c_R^{\vert D\vert} (  A^*,B) . $$ Whenever $A$ and $B$ are
Hilbert-Schmidt, then
$$\bar c_S^D(j(A), j(B))= 2i \omega^D (A,B)=
2tr(A^*B-B^*A) $$ where $tr$ is now an ordinary
trace.  }
\bp
\bp {\bf 5. Geometry on weighted bundles}
\mp $\bullet$ {\it Weighted vector bundles}
  \sp  We shall say that a  Hilbert space $H$ lies
in the class
${\cal C}{\cal H} $   whenever there is   a
compact boundaryless compact smooth Riemannian
manifold $M$, 
  a finite rank smooth vector bundle $E$ based on $M$ 
 and  $s>{dim M \over 2}$ such that
  $H=H^s(M,E)$. 
Typically, letting $G$ be a Lie group and $Lie(G)$ be its Lie algebra, 
 the Lie 
 algebra 
$H^s(M, Lie (G))$ of the Hilbert current group $H^s(M, G)$
 lies in ${\cal C} {\cal H}$.
\sp 
 Let 
 ${\cal C} {\cal E}$ be the class of Hilbert vector bundles   ${\cal E}\to X$
based on a 
  (possibly infinite dimensional) manifold $X$ 
with    fibres    modelled 
on a separable Hilbert space $H$ in the 
class ${\cal C}{\cal H} $ defined   previously   and  with  {\it
transition maps   in $PDO(M, E)$} when  $H $ is the space of sections in some
Sobolev class of a vector bundle $E$ based on $M$. \sp 
 ${\cal C}  {\cal X}   $  denotes the class 
of   infinite dimensional manifolds $X$
 with
tangent bundle $TX$   in ${\cal C} {\cal E} $. Since the transition maps are bounded, 
 they correspond to operators of order $0$ and since they are invertible, they 
are in fact elliptic operators of order zero so that they lie
 in $Ell(M, E)$. 
\mp To illustrate this setting   let us   give
examples of   manifolds, resp.   vector bundles in the class ${\cal C} {\cal X}$, resp. 
 ${\cal C}{\cal E}$.
 \mp $\bullet$ {\it Examples  }\sp
i)  Finite rank vector bundles  lie in the class ${\cal C}{\cal E}$. To see this we take 
  a manifold $M=\{\ast \}$ reduced to  a point $\ast$, the bundle $E$ to be   trivial of the type 
 $\{\ast\}\times \R^d$ (or $\{\ast\}\times \C^d$ if the bundle is complex). 
The transition functions 
belong to $Ell (\{\ast\} , E)= Gl_d (\R)$ ( or $Gl_d (\C)$ if the bundle is
complex).
\sp  ii) If $G$ denotes a Lie group and $s>{dim M \over 2}$, then the current group
 $H^s(M, G)$ is a Hilbert Lie group and can be equipped with 
a left invariant atlas $\phi_\gamma(u)(x):=  exp_{\gamma(x)}
( u(x))\quad \forall x\in M, \quad \forall \gamma \in H^s(M, G)$
where $exp_{\gamma(x)}$ is the exponential coordinate chart at point $\gamma(x)$
induced by a 
left invariant Riemannian metric  on $G$. The transition functions 
are given by multiplication operators which indeed are P.D.Os.
 \sp iii) More generally, let $N$ be a Riemannian
manifold, then the space $H^s(M, N)$ is a Hilbert
manifold with tangent space at a point
$\gamma$ given by $H^s(M,
\gamma^* TN)$. This manifold, which is modelled
on  $H^s(M, \R^n)$ where $n$ is the dimension of
$N$,  can be equipped with an atlas induced  by the
exponential map $exp^N$ on
$N$ in a similar way to the above
 description, a
local chart being of the type
$\phi_\gamma(u) (x)= exp^N_{\gamma(x)} (u(x))$.
The transition functions are locally given by
multiplication operators and hence define PDO's.
\mp $\bullet$ {\it Bundles of
operators }\sp 
 Let ${\cal E}$ be a $\Ci$
 vector bundle in the class ${\cal C}{\cal E} $
based on a manifold $B$  and let it be 
modelled on a separable Hilbert space 
$H$.
 For any $b\in B$ we shall denote by  $PDO_b({\cal
E} )$ 
 the class of operators 
  $A_b $ acting densely on the fibre ${\cal E}_b $
above $b$
  such that    for any local trivialisation 
$\phi: {\cal E}\vert_{U_b}  \to U_b\times  H$
around $b$ the operator 
$  \phi^\sharp  A(b):= \phi(b)A_b \phi(b)^{-1} $
(where $\phi(b): {\cal E}_b
\to H$ (obtained after localizing it using smooth cut-off functions) is the 
isomorphism induced by the trivialization)
 lies in $PDO(M, E)$. This definition  involves
local charts but it is in fact independent of the
  choice of local chart. Indeed,   for another local chart $(U,\psi)$ we
have 
  $ \psi^\sharp  A(b) =\psi(b) \circ \phi(b)^{-1}
\phi^\sharp A_b \phi(b) \circ
\psi(b)^{-1}$,  and since
  the transition functions $\psi(b) \circ
\phi(b)^{-1}$ are given by classical P.D.Os, 
the condition 
$\phi^\sharp A  \in PDO(M,E) $
 is independent of the choice  of $\phi $.  In a similar way, 
the notion of order defined by the order of the P.D.O in the local chart,
does not depend 
on the choice of local chart.  
   \sp We shall also denote by  $Ell_b({\cal E} )$
the class of    operators
$A_b $ acting densely on 
${\cal E}_b$
 such that     for any local trivialisation $\phi: {\cal E}/U_p  \to U_p
\times  H$ around $b$ the operator 
 $\phi: {\cal E}/U_b   \to U_b \times  H$ around
$b$ the operator, $(\phi^\sharp A)(b)   $ lies in
$Ell(M, E)$.  Here again, the definition involves
a choice of local chart but is in fact
independent  of that choice.Indeed, the principal
symbol being  multiplicative and using the
characterization of ellipticity  in terms of 
invertibility of the principal symbol, one easily
checks that the condition  $\phi^\sharp A \in
Ell(M, E) $
 is   independent of the choice of   $\phi $.  
    \sp 
 \sp This  gives rise to two  bundles $PDO( {\cal
E}):= \bigcup_{b \in B} PDO_b({\cal E} )$ and 
 $Ell({\cal E}):= \bigcup_{b \in B}
 Ell_b({\cal E}
)$  with fibre at point
$b$ given 
respectively by 
  $PDO_b( 
 {\cal E} )$    and   
 $Ell_b({\cal E} )$.  
 \sp Notice that when ${\cal E}$ is a bundle of 
finite rank,
we have $PDO({\cal E})= Hom({\cal E})$ and 
 $Ell ({\cal E})= GL({\cal E})$ since the underlying manifold $M$ 
 reduces to a point and the
 vector bundle to a vector space, namely the model space of ${\cal E}$.  
\mp 
$\bullet$ {\it Weighted traces}
\sp A weight on a  smooth  finite rank   
vector bundle in the   class ${\cal C} {\cal E}$
is a smooth section of
$Ell({\cal E})$ of operators with constant
order which  is locally 
   a positive
self-adjoint (elliptic) operator. 
  If ${\cal E}$ has finite rank, it  is simply
given
by a section of
$GL({\cal E})$ which locally 
 is a positive
self-adjoint operator   and hence
corresponds to the choice of a Riemannian
metric.
\sp If  the base manifold is 
either locally compact and a countable union of
compacts 
or if it is a paracompact Hilbert manifold
 modelled on a
 separable Hilbert space, then it has a
smooth partition of unity [L] and such a
global section can be built up patching up local
sections of 
$Ell({\cal E})$, i.e maps from an open
neighborhood in the base manifold $B$ to $Ell(M,
E)$  (this construction is similar to the
one that gives the existence of a Riemannian metric
under the same conditions).
\sp  Let
${\cal E}\in {\cal C} {\cal E}$ be modelled on
$H$. We define the weighted pseudo-trace of a
field of operators:
   $$\eqalign{ \Gamma(PDO({\cal E}))&\to  \Gamma(X, \R)  (\hbox{or}\quad \Gamma(X, \C)) \cr
A&\to tr^Q(A) \cr} \eqno(5.1)$$
 locally; in a local chart $(U, \phi)$ it
 is defined by
the $\phi^\sharp$-weighted pseudo-trace 
$tr^{\phi^\sharp Q}(\phi^\sharp A )$ of
$\phi^\sharp A$. Although the definition involves
a choice of local
 chart, because of the  
covariance
property   it is  in fact independent of this
choice. Indeed  if 
$\phi $ and $\psi $ are two local
 charts around
$x$ and  if we set 
$C := \psi   \circ  \phi^{-1}$  in  formula
(1.7),   we have:
$$   tr^{\psi^\sharp Q} (\psi^\sharp A)    = 
tr^{C\phi^\sharp Q C^{-1}} (C\phi^\sharp A C^{-1})=
tr^{\phi^\sharp Q} (\phi^\sharp A).  $$   \mp 
 {\it Remark }\sp    Looking back at definition
(1.3), strictly speaking, 
 in order to make sure the field $P_Q$ of
orthogonal projections is smooth, 
     we should assume that 
$KerQ$ has constant dimension.   
This    is an artificial difficulty   which comes  from the specific construction of the zeta
 function renormalization  since it 
 involves taking complex powers of an invertible operator 
$  (Q+P_Q)^{-z} $. However, bearing in mind that 
 the heat-kernel  renormalization only involves  the exponential $ 
e^{ -tQ}$   and the  fact that $Lim_{z\to 0}tr(A (Q+P_Q)^{-z})=Lim_{\e \to 0} 
tr(A e^{-\e Q})$ (which can be shown via the Mellin transform), we see that the jumps in the dimension of the kernel of $Q$ do not affect the
renormalized trace. From now on, we shall assume that $Q$ is invertible, otherwise one just replaces $Q$ by $Q^\prime:= Q+P_Q$ where $P_Q$ is the orthogonal projection onto the kernel of $Q$. 
 \mp $\bullet$ {\it Variations of 
 weighted traces}
\sp  {\bf  Proposition 4} {\it
\sp   Let
$({\cal E},Q)$
 be a weighted vector bundle
equipped with a connection 
$\nabla$ and let $\alpha$, $\beta$ 
 be two  $ PDO({\cal E} ) $
valued   forms on the base 
manifold of ${\cal E}$. 
\sp  1)    $$tr^Q([\alpha, \beta])= -{1\over ord Q}
res ([log Q, \alpha]
\beta)\eqno(5.2) $$
\sp  2)   Provided both $[\nabla, log Q]$  and $
[\nabla, \alpha]$ are sections of $PDO({\cal E}),
$ we have
$$   [\nabla, tr^Q] (\alpha):=  dtr^Q(\alpha) -=
tr^Q([\nabla,\alpha])=-{1\over ord Q} res(\alpha
\cdot [\nabla,    log Q]). \eqno(5.3)$$
Here $ord Q$ denotes the  (constant) order of the field of elliptic
operators $Q$ 
.}
  \sp   {\it Proof} \sp   (5.2 ) follows from
Proposition 1 (compare with (2.2)) 
\sp As in the proof of Proposition 1, to prove 
 (5.3) we use   the fundamental property  of the 
canonical trace
$TR$.    Let us   first consider  a smooth
(for the natural topology
on classical P.D.Os induced by the topology of 
uniform convergence in all derivatives  on the classical symbols) one
parameter family of operators
$A_t\in PDO(M, E)$ with  constant order $a$ and a smooth one parameter
 family $Q_t\in Ell^+_{ord
>0}(M, E)$ (for the natural Fr\'echet topology on  $PDO(M, E)$)
with constant order $q$, $t $ varying in 
$ ]0, 1[$.  Applying formula  (1.3) to
$A_z:= { A_0 (Q_t^{-z}- Q_0^{-z})\over z}$ (in
which case $\alpha(z)= a-q z$)  and then going to
the limit when $t\to 0$  yields:  
$$   {d\over dt}_{/t=0} tr^{Q_t}(A_0) =  -{1\over q}
 res(A_0  {d\over dt}_{/t=0} log Q_t).
$$ let us now  also consider a $1$-parameter smooth
family $(A_t)$ of P.D.Os of constant order, then:
$$   {d\over dt}_{/t=0} tr^{Q_t}(A_t) = tr^{Q_0}({d\over
dt}_{/t=0}A_t)-{1\over ord Q_0}
 res(A_0  {d\over dt}_{/t=0} log Q_t).
$$
Similarly, given any local trivialization around a point $b_0$ in the base
manifold $B$ of 
${\cal E}$, we have:
$$   dtr^Q(\alpha) = tr^Q(d\alpha )-{1\over ord Q} res(\alpha \cdot d   
log Q  ). \eqno(5.4)$$
Let us write $\nabla= d+ \theta$ in this local trivialization.  Since by
assumption 
$[\nabla, \alpha]$ is a section of $ PDO({\cal E})$, we have that locally 
$d \alpha + [\theta, \alpha] \in PDO(M, E)$  and hence, since $d\alpha\in
PDO(M, E)$ as the
 differential of a P.D.O, we conclude that 
$ [\theta, \alpha]$ also lies in $PDO(M, E)$. 
Applying the first part of the lemma yields:
$$tr^Q([\theta, \alpha])= -{1\over ord Q} res (
[log Q, \theta]
\alpha).\eqno(5.5)$$ 
Combining (5.4) and (5.5) we find:
$$ \eqalign{   dtr^Q(\alpha) &= tr^Q(d\alpha )-{1\over ord Q} res(\alpha
\cdot d   
log Q  )\cr
&= tr^Q([\nabla, \alpha] )- tr^Q ([\theta, \alpha])-{1\over ord Q}
res(\alpha \cdot d   
log Q  )\cr
&= tr^Q([\nabla, \alpha] )+ {1\over ord Q} res([log Q, \theta]\alpha
)-{1\over ord Q} res(\alpha 
\cdot d   
log Q  )\cr 
&= tr^Q([\nabla, \alpha] ) -{1\over ord Q} res(\alpha 
  [\nabla,    
log Q] ).\cr}\hfill \bullet$$ 
\mp $\bullet$ {\it Ricci curvature on weighted
manifolds }\sp Let $X\in {\cal C} {\cal X}$ be
equipped with  a Levi-Civita connection $\nabla^X$
such that for any tangent vector fields  $U, V $ 
the map:
$$ R^X(U, V): W\to \Omega^s(W, U) V \quad \forall  W\in \Gamma(TX)$$ 
 is a section of  $ PDO(TX)$ where $\Omega^X$ denotes the curvature of $\nabla^X$.  
Then, given a weight $Q$ on $X$ we can define the $Q$-{\it weighted Ricci curvature}:
$$r_1^{X, Q} := tr^Q (R^X ).\eqno (5.6)$$
Since $tr^Q$ coincides with the ordinary trace
 on trace-class operators if $X$ is finite
dimensional, it then coincides with the ordinary
Ricci curvature. 
 \mp\mp   $\bullet$ {\it
 Weighted first Chern form  on
K\"ahler manifolds }\sp Let us now assume  
$X\in {\cal C} {\cal X}$ is K\"ahler and let 
 $\nabla^X$  be the a K\"ahler  connection, let 
$Q$ be a weight on $X$. In a similar way, provided
that for any holomorphic   tangent vector  field 
$U$ and for any antiholomorphic tangent vector
field 
$ \bar 
V $ the map $\Omega^X(U, \bar V)$ is a section   of  $ PDO(TX)$, we define 
the $Q$-{\it weighted first Chern form}: 
 $$r_1^{X, Q}(U, \bar V) := 
tr^Q (\Omega^X (U, \bar V)) \eqno (5.7)$$ 
where $\Omega^X$ denotes the curvature of $\nabla^X$. If $X$ is finite dimensional, it coincides with 
the ordinary first Chern form.\sp 
 Unlike the first Chern form on a   finite dimensional manifold,
it is not closed in general.  This follows from the following proposition: 
   \mp   {\bf Corollary 4}\sp {\it Let
$(X,\nabla^X)$ be
 a K\"ahler manifold eqXipped with a
 weight $Q$    and let
$ord Q$ be the (constant) order of $Q$. Then,
provided  $[\nabla^X, log Q]$ is a $PDO(TX)$
valued one form and provided the curvature
$\Omega^X$ is a 
$PDO(TX)$ valued two form, we have 
 $$ dr_1^Q =    -{1\over ord Q} res  
 ([\nabla^{X}, log Q]\Omega^X). 
\eqno (5.8) $$} 
\sp {\bf Proof:} Applying Proposition 4 to
$\alpha:=
\Omega^X$, we find:
$$\eqalign{ dtr^Q(\Omega^X)&= tr^Q([\nabla^X,
 \Omega^X]) 
 - 
{1\over ord Q} res (\Omega^X [\nabla^X, 
 log Q ])\cr
&= - 
{1\over ord Q} res (\Omega^X [\nabla^X,  log Q ])\cr}$$
where we have used the Bianchi identity 
  $[\nabla^X, \Omega^X]=0$.
 \hfill $\bullet$
\mp 
$\bullet$ {\it Remark  } \sp From Kuiper's
results [Ku] on the contractibility of  the
unitary group of a separable Hilbert space, we
know that the orthonormal frame bundle 
$O({\cal E})$ of a hermitian  vector bundle
 with fibres modelled on a Hilbert space, is
topologically trivial. This is the case for the
class of manifolds we are investigating, and 
from this topological triviality one might expect
that whenever the first Chern class is closed,
teh corresponding  characteristic class should
vanish. This is not the case as we shall see
shortly in the case of current groups, but in our
setting  the non triviality seems to come  from
the fact that we restict ourselves to  P.D.Os.
Indeed, the holonomy bundle is a reduction of the
frame bundle and for a manifold
$X$  in the class ${\cal C} {\cal X}$ (recall that
we only allowed transition maps which were 
P.D.Os) modelled on the
 Hilbert space $H:= H^s(M, E)$ 
for some $s>{dim M \over 2}$, 
  the  structure group of the frame bundle is 
  $$  GL^{Ell}(H) := GL(H)\cap  PDO(M, E)= 
GL(H)\cap Ell(M, E) = Ell_0^*(M, E)\eqno (5.9)$$
where $Ell_0^*(M, E)$ denotes the group of invertible zero order elliptic 
operators acting on sections of
 $E$. But it is
a well-known fact that the pathwise connected component of identity of 
$ Ell_0^*(M, E)  $ has a non trivial
fundamental group isomorphic to   $K_0(S^*
M)$ where $S^*M$ is the unit sphere in the
cotangent bundle of $M$ (see Appendix B).  
\mp  
   \bp \noi {\bf 6. Weighted Lie  groups}
  \mp In this section, we apply the constructions
 and results of  the previous section 
 to the case of weighted   Lie groups, thus preparing for the next section where we will specialize to
   current groups. Here ${\cal G}$ is an
infinite dimensional Hilbert Lie group  in the class ${\cal C} {\cal X}$ with 
Lie algebra $Lie ({\cal G})= H^s(M, E)$
(for some $s>{dim M \over 2}$ and  some hermitian vector bundle $E$ based on some manifold $M$). 
\mp  $\bullet $ {\it Left invariant weights
}\sp     A natural  weight on ${\cal G}$ is given
by  a {\it left invariant field} of operators
$$Q(\gamma):=  {L_\gamma}_* Q_0
{L_\gamma}_*^{-1}\quad \forall
 \gamma \in {\cal G}\eqno(6.1)$$ 
 where $Q_0 \in Ell_{ord>0}^{ +}(M, E)$ is any weight
 on the Lie algebra $ Lie ({\cal G})$   and where 
${L_\gamma}$ denotes left multiplication.  
\sp 
 As a consequence of Proposition 4 we have: 
\smallskip \noindent 
\noindent {\bf Corollary 5}\smallskip \noindent
{\it Let $\nabla$ be a left invariant connection on
${\cal G}$ induced by a left invariant one form
$\theta_0\in Lie({\cal G}) \otimes Hom(Lie({\cal
G}))$. Then, given  
 a left invariant 
$p$-form $\omega$ on
 ${\cal
G}$ which we identify with $\omega_0\in
\Lambda^p(Lie {\cal G}) \otimes PDO(M, E)$, we
have: 
$$[\nabla, tr^Q] (\omega):= d tr^Q(\omega)-
tr^Q([\nabla, \omega])= -tr^{Q_0} ([\theta_0,
\omega_0])
 = {1\over ord Q_0} res( [log Q_0,
\theta_0] \omega_0) $$
where the bracket is an operator bracket so that 
 $[\theta_0, \omega_0](X, X_1, \cdots, X_p)=
[\theta_0(X), \omega_0(X_1, \cdots, X_p)].$ }
\smallskip \noindent 
{\it Proof:} Recall that for a left invariant
$p$-form $\alpha$ on ${\cal G}$ we have:
$$d\alpha(X_1, \cdots, X_{p+1})= \sum_{i<j}
(-1)^{i+j} \alpha([X_i, X_j], X_1, \cdots, \hat
X_i, \cdots, \hat X_j, \cdots, X_{p+1})$$
 whereby the $X_i's$ are left invariant vector
fields. Applying this to $\alpha= tr^Q(\omega)$
and $\alpha= \omega$ and then taking the
$Q$-weighted trace yields $dtr^Q(\omega)= tr^Q(d
\omega)$. Hence
$$\eqalign{ [\nabla, tr^Q] (\omega)&:= d
 tr^Q(\omega)- tr^Q([\nabla, \omega])\cr
&= 
-tr^Q ([\theta, \omega])\cr
&= -tr^{Q_0} ([\theta_0, \omega_0])\cr
&= {1\over ord Q_0} res( [log Q_0,
\theta_0] \omega_0).\cr}$$ 
\mp 
$\bullet$ {\it Diagonal weights }\sp  Let us now
consider the particular case when  $E$ is  a
trivial bundle $E= M\times V$ where $V$ is a 
finite dimensional  hermitian vector space. We  
can restrict ourselves to  (left invariant) {\it
diagonal weights},
 meaning by this that: 
$$Q_0:= \bar Q_0\otimes 1_V \eqno(6.2)$$
where $1_V$ is the identity opertor on $V$ so
that (6.1) reads: 
$$Q(\gamma)= L_{\gamma_*}\left(\bar  Q_0\otimes
1_V\right) L_{\gamma_*}^{-1}.\eqno(6.3)$$  In that
case, for a left invariant
$p$-form induced by $\omega_0$ as in the above
corollary we have:
$$\eqalign{ tr^Q(\omega)&= tr^{Q_0} (\omega_0)\cr
&= Lim_{ z\to 0}tr( \omega_0(
Q_0+P_{Q_0})^{-z})\cr &= Lim_{ z\to 0}tr( \omega_0(
\bar Q_0+P_{\bar Q_0})^{-z}\otimes 1)\cr
&= Lim_{ z\to 0}tr( (
\bar Q_0+P_{\bar Q_0})^{-z}tr_V(\omega_0))\cr
  &= tr^{\bar Q_0} (tr_V(\omega_0))
\cr}\eqno (6.4)$$ where
$tr_V$ denotes the  
 finite dimensional trace  on the Lie algebra $V$.  
Hence,  provided $tr_{V}(\omega_0)$ is itself a
trace-class operator (i.e a P.D.O of order
$<-{dim M\over 2}$) 
then $$tr^Q(\omega) =
tr(tr_{V}(\omega_0))\eqno(6.5)$$ (where
$tr$ is an ordinary trace) is independent of $\bar
Q_0$.  This is the type of  situation we shall be
working  with in these notes.  Since there is no
weight dependence in that case, we can choose any
operator $\bar Q_0 $ and in particular  
 the Laplace operator $\Delta$ on the Riemannian
 manifold  
$M$.    \bp  {\bf 7. The case of current groups }
\mp   We now specialize to the current group 
$H^s(M, G)$ ($s>{dim M \over 2}$) of $H^s$ maps
from the  compact Riemannian manifold $M$ to a
semi-simple  Lie group $G$ of compact type
  (which ensures that the Killing
form is non degenerate and that the adjoint
representation $ad$ on the Lie algebra is
antisymmetric for this bilinear form). $H^s(M,
G)$  is a Hilbert Lie group with     Lie algebra  
$H^s(M, Lie (G))$ where $Lie (G)$ denotes the 
Lie algebra of $G$ so that the current group $H^s(M, G)$ is a manifold in the class 
  ${\cal C} {\cal X}$, the underlying finite rank vector bundle $E$ being trivial 
since $E= M\times Lie(G)$ and we     
equip $H^s(M, G)$ with the 
{\it left invariant and
diagonal weight} $Q$ (see (6.2) with $V=Lie(G)$): 
 $$Q(\gamma):= L_\gamma ( \bar Q \otimes 1_{Lie(G)})L_\gamma^{-1}= L_\gamma 
  \left(\Delta \otimes 1_{Lie(G)}\right)
 L_\gamma^{-1}\quad \forall \gamma \in H^s(M, G).
\eqno(7.1)$$ 
\mp 
$\bullet$ {\it A left invariant metric }\sp   The
operator $Q_0:= \Delta \otimes 1_{Lie(G)} 
 $    is an elliptic operator of order $2 $ acting
densely on $ H^s(M, Lie(G))$. It is 
   positive for the scalar product 
$$\langle \cdot, \cdot \rangle_0:=  
\int_M dvol(x) (\cdot, \cdot)_{Lie(G)}\eqno(7.2)$$
where $(\cdot, \cdot)_{Lie(G)}$ is a scalar product on $Lie(G)$ given by minus
 the Killing form.
$H^s(M, G)$ is equipped with a  left invariant
metric   defined in terms of the following scalar
product on $H^s(M, Lie(G)) $: 
 $$\langle \cdot, \cdot \rangle_0^s:=   
\langle {(Q_0+P_{Q_0})}^{s\over 2} \cdot,
 {(Q_0+P_{Q_0})}^{s\over 2}\cdot
\rangle_0.\eqno(7.3)$$
\mp $\bullet$ {\it The Levi-Civita connection }\sp
The corresponding the Levi-Civita connection is  
a left invariant connection
     $\nabla^s   = d+ \theta^s$
where $\theta^s$ is a left invariant
 $Hom(H^s(M, Lie(G))$ valued one form on
$H^s(M,G)$ given  by a map $\theta_0^s: Lie({\cal
G})\to Hom( Lie{\cal G})$  (see [F1] formula (1.9)
up to a sign mistake):
$$\theta_0^s(U)= {1\over 2} \left(ad_U+
{(Q_0+P_{Q_0})}^{-s} ad_U 
  (Q_0+P_{Q_0})^s - 
(Q_0+P_{Q_0})^{-s} ad_{
(Q_0+P_{Q_0})^sU}\right)\eqno (7.4)$$ where $U$ is
 an element of the Lie algebra
 $H^s(M, Lie(G))$.  \mp 
 $\bullet$ {\it The   curvature }\sp The curvature
$\Omega^s
$ is a left invariant two form given by  an 
element
$\Omega_0^s\in
\Lambda^2(H^s(M, Lie(G))\otimes Hom(H^s(M,
Lie(G))$:
$$\Omega_0^s(U, V)=  \theta_0^s\wedge
\theta_0^s(U, V)= [\theta^s(U), \theta^s(V)]-
\theta_0^s([U, V])\quad U, V\in
H^s(M, Lie(G)).\eqno(7.5)$$  
\mp {\it Warning} \sp {\it We shall henceforth
identify left invariant forms with  forms on
the Lie algebra $H^s(M, Lie(G))$ so that for
 left  invariant fields $X, Y$ induced
respectively by   elements
$X_0, Y_0\in H^s(M, Lie(G))$ we have  
$\theta^s(X )= \theta^s_0(X_0 )$
and $\Omega^s(X, Y)= \Omega^s_0(X_0, Y_0)$. } 
\sp   D.Freed shows ([F1] prop. 1.14) that  for smooth
$X, Y\in H^s(M, Lie(G))$  
 the  map $$R^s(X,Y): Z\to \Omega^s(Z, X) Y$$ is a pseudo-differential operator 
(with Sobolev coefficients) 
of order 
$max(-1, -2s)$. Freed in fact proves the result
 for smooth $X$ and $Y$ but his
 proof easily extends to the case when $X, Y\in H^s(M, Lie(G))$ since it involves a counting 
of the order of the operator $R^s(X,Y)$
 which is independent of the degree of regularity
of its coefficients. 
 Notice that for $s\in {1\over 2} \N-\{0\}$, we have
 $max(-1,-2s)=-1$. \sp  We further quote from
[F1]: \sp \noi  {\bf Lemma 5}\sp\noi  {\it  For
$X, Y\in H^s(M, Lie(G))$ the operator 
$tr_{Lie(G)} R^s(X, Y)$ is a classical
pseudo-differential operator of order $-2q  $ where
 $q= min (1, 2s)$ and where as before 
$tr_{Lie(G)}$ denotes the finite dimensional trace
on $Lie(G)$. }
\mp \noi {\bf Proof}: We refer the reader to    [F1] Proposition 1.18  where  this result  was proven 
in   the case when   $X, Y\in \Ci(M, Lie(G))$.  
\hfill $\bullet$
\mp $\bullet$ {\it The Ricci  curvature }\sp 
Following the general procedure described
in section 5, let  us define a weighted 
 Ricci curvature on the   Lie group $ H^s(M, G) $  and compare it to that of [F1].\sp \noi 
For any diagonal weight $Q$, as a consequence of
(6.5) we have:  
 $$ Ricc^{s,Q}(X,Y):= tr^Q(R^s(X, Y))=
tr^{\bar Q_0}\left(tr_{Lie ( G)}(R^s(X,
Y))\right). \eqno( 7.6)$$
 Hence provided $min(1,2s)  >{dim M\over 2}$, using 
  Lemma 5 we can write: $$ Ricc^{s,Q}(X,Y)=     
tr( tr_{Lie(G)} (R^s(  X, Y  ) )) \eqno(7.7) $$
wher $tr$ is now an ordinary trace, since in that
 case, 
 $tr_{Lie(G)} (R^s(  X, Y  ) )$ is trace-class.  This holds in particular  for any
 {\it loop group} 
$H^s(S^1, G)$ with $s>{1\over 2}$. 
 \sp It is easy to check  that in that case and 
for any left invariant weight $Q$  the Wodzicki
residue 
 $res(  R^s(X, Y) )$ of the operator $R^s(X, Y)$
 vanishes and that the $Q$-weighted Ricci
curvature is given by an ordinary limit: 
$$Ricc^{s, Q} (X, Y)= lim_{z\to 0} (TR( R^s(X, Y)
Q^{-z}))\eqno(7.8)$$ however this limit depends on
the choice of $Q$. 
\mp 
$\bullet$ {\it The case when
$s$ is an integer and $M$ is odd dimensional }\sp
There are some cases for which the
$Q$-dependence vanishes; we shall need a
preliminary lemma to single them out. 
 \sp \noi 
{\bf Lemma 6}  \sp\noi  {\it  For two left
invariant vector fields  $X, Y$ on $  H^s(M,
 G )$ and 
$s\in
\N$  the operator
 $R^s(X, Y)$     lies in the odd-class.    }
\sp \noi {\it Proof} \sp \noi 
From the expression of the Levi-Civita curvature
 (7.5)
one sees  that 
the operator  $R^s(X,Y)$ is built up from
compositions and linear combinations of  operators
in the odd class provided $Q_0$ is in the odd
class, since $ad_X$ and $adX$   both lie in the odd
class since they are built up from multiplication 
 operators. Since for any integer $s$, the operator  $Q_0^s$ does lie   in the odd 
class  and since the product of  odd-class operators lies in the odd class,  the result follows. 
 \bp \noi 
 \mp {\bf Proposition 5}\sp \noi  
{\it  When the underlying manifold is odd
 dimensional 
and $s$ is an integer,   the (left invariant) 
weighted  Ricci curvature on  $H^s(M, G)$
equipped with the left invariant connection (7.4) 
is independent of the choice of the  weight
$Q$ among operators in the odd class
 and we have for any two left invariant vector
fields $X$ and $ Y$ on $H^s(M,  G )$:
$$Ricc^{s, Q }(X,Y) = TR_{odd} (R^s(X,
Y)).\eqno(7.9)$$  }
 \mp {\bf Proof }: 
\sp\noi  
As we saw in section 1 (formula (1.6)),  the
dependence on the choice of $Q$ is   measured  in
terms of a residue since for two weights $Q_1$ and
$Q_2$ we  have:$$ tr^{Q_1}(R^s(  X, Y)) - 
tr^{Q_2}(R^s(  X, Y))= -q^{-1} \left(res ( R^s( 
X, Y) (log Q_1-log Q_2))
\right) $$ where $q$ is the common order of $Q_1$
and $Q_2$.  If both $Q_1$ and $Q_2$  lie  in the
odd-class so does $  log Q_1-log Q_2 $  (see Proof
of Proposition 4.1 in [KV1]), and  since $R^s(X,
Y)$ also  lies in the odd class,  the operator
$R^s(  X, Y) (log Q_1-log Q_2)) $ lies in the
odd-class. 
 Thus in odd dimensions its Wodzicki  residue vanishes. \hfill $\bullet$
\mp \noi   There is a priori no reason for a
similar property to hold for $s={1\over 2}$, which
is the  case we shall be focusing on in the
following.   The $Q$-weighted Ricci curvature will
in general depend on the choice of $Q$. \bp
  {\bf 8.  The case of  the
based loop group $H^{1\over 2}_e(S^1,G)$ }\mp We
now specialize to the case when $M=S^1$ and
$s={1\over 2}$, applying the results of the
previous section to  a based loop group. Notice
that for
$s={1\over 2}$, $H^s$ maps are not even continuous
in  that case, and in practice, one works  with 
$H^{{1\over 2}+\e}, \e>0$ in order to have continuous objects.\sp The space 
 $H_e^s(S^1, G)\subset H^s(S^1, G)$ of  $G$ valued  loops  with value $e_G$ at a given point,
 where $e_G$ is the identity element in $G$. It is a Hilbert manifold with Lie algebra
the based loop algebra  
$H_0^s(S^1, Lie(G))$ of maps in the loop algebra which all coincide at  point $0 $.  
$H_e^s(S^1, G)$
 is equipped with an almost complex structure, for
which the metric 
 is in fact K\"ahler when  $s={1\over 2}$ [F1]. We
 equip it  as before with a left invariant
diagonal weight  
$Q$, the operator $Q$ being the same one used to
define the connection $\theta^s$.
 \mp  $\bullet $ {\it An almost complex 
structure    on $H^{1\over 2}_e(S^1, G)$ }\sp 
   We introduce a  Dirac  operator 
$\bar D_0=z{ d\over dz} =
-i{d\over d t}
 $ (with $z=e^{it}$) acting  on
   $\Ci(S^1,\C)$ and we set $D_0:= \bar
D_0\otimes 1_{Lie(G)}= z{ d\over dz}\otimes
1_{Lie(G)}$, $D$ denoting the left invariant
field of operators generated by $D_0$.  We shall
choose
$Q_0:=  D_0^2= \Delta \otimes 1_{Lie (G)}$
where $\Delta$ is the Laplacian on functions,  to
 define the left invariant
diagonal weight $Q$.   $D_0$  is injective when
restricted to based loops and the sign 
$\e(D_0):=  D_0 \vert D_0\vert^{-1}$ of the Dirac
operator is a pseudo-differential operator of 
order $0 $ which yields a conjugation on 
$H_0^{1\over 2}(S^1,  Lie(G))$ since
$\e(D_0)^2=1$.
 We have the splitting:
$$H:=H_0^{1\over 2}(S^1,  Lie(G))=H_+  \oplus H_-
\eqno(8.1)$$ 
 where  $$H_+ :=Ker(\e(D_0)-1)=\pi_+ 
\left(H_0^{1\over 2}(S^1, Lie(G))\right)$$   
  and $$H_- :=Ker(\e(D_0)+1)=\pi_-  
\left(H_0^{1\over 2}(S^1, Lie(G))\right).$$ Here
$\pi_+  $ and $\pi_+  $ are the orthogonal projections 
  w.r. to the scalar product  
 $\langle\cdot, \cdot \rangle_0^{1\over 2}$
 defined in (7.3)
  $\pi_+:= {\e(D_0) +
1\over 2}$,$\pi_-:={ -\e(D_0)  + 1\over 2}$.\sp
 The set $\{u_n(t):=e^{ i    n t}, n\in \Z\} $
 yields a C.O.N.S of eigenvectors of $D_0$
 corresponding to the set of eigenvalues 
$\{\l_n:=      n, n\in \Z\}$. Then the set 
$\{u_n^+(t):=e^{ i     nt}, n\in \N\} $ 
spans $H_+$ and the set 
$\{u_n^-(t):=e^{ -i   nt}, n\in \N\} $ 
spans $H_-$. 
The map $J_0$ 
defined by $J_0u_n^+:=-i u_n^+$, $J_0u_n^-:= 
i u_n^-$ or equivalently by $J_0:=i\e(D_0)$  obeys
the relation   $J_0^2 =-1$ and 
  yields
a natural  almost complex structure on $H$ for which the $(1,0)$ part  $H^{1,0}:= Ker(J+i) $  
 coincides with 
$H_+$, the $(0, 1)$ part  $H^{1,0}:= Ker(J-i)=H_+$  with $H_-$. 
  By left invariance of the metric 
on $H_0^{1\over 2}(S^1, Lie(G))$, this gives rise
to  a left invariant almost complex structure
 $J(\gamma):= L_\gamma J_0 L_\gamma^{-1}$ 
on the Lie group 
$H_e^{1\over 2}(S^1,  G) $. \mp 
 $\bullet $ {\it The K\"ahler connection     on
$H^{1\over 2}_e(S^1, G)$ }\smallskip \noindent  
$H^{1\over 2}_e(S^1, G)$ is equipped with a left
invariant symplectic   (hence closed see [Pr])
form: 
$$\omega(X, Y):=  \int_0^1\langle
X_0^\prime, Y_0\rangle d t  $$
where $X, Y$ are two left invariant vector
fields generated by $X_0, Y_0$ and where the
"prime" denotes derivation with respect to $t$. One
can check that the  invariant bilinear form given
by
$B( X, Y):= \Omega(X, JY)$ yields  back the
$H^{1\over 2} $ Sobolev metric  given by (7.3).  The 
associated left invariant hermitian form reads: 
 $$\langle\langle \cdot  , \cdot \rangle \rangle_0^{1\over 2}:= 
\langle \cdot , \e(D_0 )\quad \cdot 
\rangle_0^{1\over 2}= \langle D_0\cdot,  \cdot
\rangle_{L^2}. \eqno (8.2)$$ where we have used the
fact that $\e(D_0)_0 Q_0^{1\over 2}=\e(D_0) \vert
D_0\vert=  D_0$ and where $\langle \cdot, \cdot,
\rangle_{L^2}$ denotes the $L^2$ hermitian product.
 Hence $S$ is  a K\"ahler form  on the
  manifold $H^{1\over 2} (S^1, G)$
equipped with the $ H^{1\over 2} $ left invariant
  Riemannian metric given by  (7.3) and the left
invariant complex structure $J$ defined above.
As a consequence, the 
 Levi-Civita connection $\nabla^{1\over 2}$
 on $H_e^{1\over 2}(S^1,  G) $ is K\"ahlerian,
 meaning by this that 
$[\nabla^{1\over 2}, J]=0$.  
\sp  Since  $\nabla^{1\over 2}$ commutes with 
$J$, $\theta^{1\over 2}$ defined in (7.4) (with
$s={1\over 2}$ )
 stablizes each of the  spaces $H_+$ and $H_-$.
Its restriction to $H_+$: 
$$\phi_0:=
\left[\theta_0^{1\over2}\right]_{++}\eqno(8.3) $$
 defines a left invariant  the (complex) K\"ahler
connection $\phi$. 
In the next lemma we    use  the 
  splitting $H= H_+\oplus H_-$ to write any  
operator
$A\in Hom (H_0^{1\over 2} (S^1, Lie(G))$ as a
matrix:   
$$A= 
\left[ \matrix {  A_{++}   & A_{+-}\cr
                              A_{-+} 
 & A_{--}\cr} \right].$$
 \sp Following [F1] we  introduce the T\"oplitz operators $T_X:= (ad_X)_{++}, \quad X\in H$. 
 \mp {\bf Lemma 7}
\sp {\it For $U\in H_+$, $$\phi(U)= D^{-1}
T_U  D .$$ For $\bar V\in H_-$, 
$$\phi(\bar V)= T_{\bar V}.$$}  \sp 
{\bf Proof:} 
\sp  1) For  $U\in H_+$, we have $\vert D_0\vert U
 =  D_0 U  . $ Hence, setting $U= \alpha \otimes
a$, $V:= \beta \otimes b$, $\alpha, \beta \in
H^{1\over 2} (S^1, \C)$, $a, b \in Lie(G)$, we
have: 
$$\eqalign{ \vert D_0\vert^{-1} ad_{\vert D_0\vert U} V&= \vert D_0\vert^{-1} (D_0 \alpha ) \beta
\otimes ad_a b\cr
&= \vert D_0 \vert^{-1} D_0 (\alpha \beta) \otimes
ad_a b- \vert D_0 \vert^{-1} \alpha D_0 \beta
\otimes ad_a b\cr
&= \vert D_0\vert^{-1} D_0 ad_U- \vert D_0
\vert^{-1} ad_U D_0\cr
&= \e(D_0) ad_U- \vert D_0\vert^{-1} ad_U
D_0.\cr}$$ Inserting this into  (7.4) yields: 
$$ \theta^{1\over 2} (U) = {1\over 2}
(ad_U+
\vert D_0\vert ^{-1} ad_U \vert D_0\vert+ \vert
D_0\vert^{-1} ad_U D_0- \e(D_0) ad_U)$$
 which, when restricted to $H_+$ and using
(8.3) leads to: 
$$\eqalign{ \phi_0(U)&= \pi_-{ad_U}_{\vert H_+}+
\vert D_0\vert^{-1} ad_U {D_0}_{\vert H_+}\cr
&=  D_0^{-1} T_U   D_0  \cr}$$
since $ad_U$ stablizes $H_+$ for $U\in H_+$. 
\sp 2)  For $\bar V\in H_-$, $\vert D_0 \vert
\bar V=-D_0 \bar V$ and hence in a similar manner
as above we have: 
$$\vert D_0\vert^{-1} ad_{\vert D_0\vert \bar V}=
\vert D_0\vert^{-1} ad_{\bar V} D_0- \e(D_0)
ad_{\bar V} $$ which yields in turn: 
$$\theta^{1\over 2} (\bar V)= {1\over 2}
(ad_{\bar V}+
\e(D) ad_{\bar V}- \vert D\vert^{-1} ad_{\bar V} D+
\vert D\vert^{-1} ad_{\bar V} \vert D\vert). $$
When restricted to $H_+$ and using (8.3), this 
reads: 
$$\phi(\bar V)= \pi_+ {ad_{\bar V}}_{\vert_{H_+}}=
 T_{\bar V} .
$$
 \bp 
 {\bf 9. The first Chern form on $H^{1\over
2}_e(S^1, G)$ and the cohomology of $PDO(S^1,
Lie(G))$ }\mp \noi 
 $\bullet $ {\it The first Chern form   on
$H^{1\over 2}_e(S^1, G)$ } \sp 
  \sp The  curvature (which also stablizes 
$H_+$ ) is given as in (7.5) by:
$$\Omega(X, \bar Y):= [\phi(X), \phi(\bar Y)]- 
\phi([X, \bar Y])=
 \left( \Omega^{1\over 2}(X, \bar Y)\right)^{1,0}.
 \eqno (9.1)$$ 
 From  [F1] (see Remarks above Theorem  2.20 ), we also know that  
 $ \Omega  (X, \bar   Y)$  
 is a classical pseudo-differential operator  
 of order  $ -1 $ and hence that it is not
 trace-class.  
 However 
 $ tr_{Lie(G)}\Omega  (X, \bar   Y)$ being of
 order $-2$ (see [F1] Remarks above Theorem 
 2.20), is trace-class. 
 We now define   the weighted first Chern form on $H^{1\over 2}_e(S^1, G)$ according to 
   (5.7) where the weighted 
K\"ahler manifold here is 
$(H^{1\over 2}_e(S^1,G),Q) $:
$$ r_1^{  Q}   := tr^Q(\Omega  )$$
Here the complex curvature is given by (9.1)
and 
   the trace is taken w.r. to the hermitian
metric. 
If the weight is diagonal then $$ r_1^{  Q}=
tr(tr_{Lie G} \Omega  ) . $$ This    follows from
(6.5) and  shows our definition of weighted first
Chern form  coincides (with a good choice of the
weight) with the "two step trace"  
  used by Freed
 [F1] to make sense of  the trace of the curvature
(see Theorem 2.20 in [F1]).    \sp  In fact only
the
$(1,0)$ part 
$Q^{1,0}:=Q_{++}$ of 
$Q$ comes into play in the 
expression $tr^Q(\Omega)$ since $\Omega$ is the $(1,0)$ part of $\Omega^{1\over 2}$.
 Before we give an expression of the first
Chern form, we need wome preliminary results. 
\mp {\bf Lemma 8} \sp {\it For the diagonal
weight 
$Q$   on $H^{1\over 2} (S^1,
G)$ chosen as in (7.1), we have 
$$tr^Q(\phi(Z))=0 \quad \forall Z\in H^{1\over
2} (S^1, Lie(G)).\eqno(9.2)$$}
\sp {\bf Proof:} It is sufficient to establish
the formula for $Z$ of the type $\gamma \otimes
c\in H^s(M, \R) \otimes Lie(G)$. If
$Z=
\gamma
\otimes c\in H_-$, then $ \phi(Z)= ({ad_Z})_{++}$
and hence (compare with (6.5))
$$\eqalign{ tr^Q(\phi(Z))&= tr^Q(({ad_Z})_{++}) \cr
&= Lim_{z\to 0}  tr\left(( {ad_Z})_{++}\bar Q_0^{-z}\otimes 1_{Lie(G)})\right)\cr
 &= Lim_{z\to 0}  tr\left(  ({M_\gamma 
\bar Q_0^{-z} \otimes ad_c  })_{++}\right)\cr
&=     tr^{\bar Q_0} ({M_\gamma  })_{++}tr_{Lie G}
(ad_c) \cr 
&= 0\cr}$$
since $ad_c$ is antisymmetric.
$M_\gamma$ denotes the multiplication operator by
$\gamma$ and $Lim$ the renormalized limit.  
  Here we have used
the fact that $Q_0$ stablizes $H_+$.
\sp If $Z\in H_+$, then $\phi(Z)= D^{-1} T_Z D =
D^{-1} {(ad_Z)}_{++} D$  and 
$$tr^Q( D^{-1} {(ad_Z)}_{++} D)=
tr^Q( {(ad_Z)}_{++})=0.$$
 \mp   {\bf Proposition 6}\sp {\it   Let
$  H^{1\over 2}_e(S^1,G)$   be equipped with a  left
invariant diagonal weight $Q$. 
  The $Q$-weighted first Chern form is the pull-back by $\phi $ of the $Q_0$ Radul
 cocycle $c_R^{Q_0}$ on $PDO(S^1, Lie(G))$ and
co\"incides with Freed's conditioned first Chern
form [F1].  For any
$X, Y\in H^{1\over 2}_0(S^1, Lie(G))$ induced by
$X_0\in H_+$ and $\bar Y_0\in H_-$ we have: 
$$    r_1^{  Q} (X, \bar Y)  
   =  c_R^{Q_0}(\phi(X),\phi(\bar Y)) 
 \eqno
(9.3)$$
  where $tr$ denotes  the ordinary trace and
$\l^D$ was defined in 
formula (4.1).}
\mp {\bf Proof:}     The first term in  the
 expression (9.1) of the curvature is a   bracket
 of P.D.Os so that taking  the weighted
pseudo-trace
   yields $tr^{Q_0} [\phi(X), \phi(\bar
Y)]=c_R^{Q_0}( \phi(X), \phi(\bar Y))$.
 The trace of the second term in (9.1) vanishes
by    lemma 7  and   we find  
$$r_1^{  Q}(X,
\bar Y)=  tr^{Q }(\Omega (X, \bar Y)) =tr^{Q_0} ([
\phi(X), \phi(\bar Y)]= c_R^{Q_0}( \phi(X),
\phi(\bar Y))$$  as announced.  \sp  \hfill
$\bullet$
\mp 
 Before we go on to the next proposition, let us 
recall various ways of  expressing the Killing
form on the Lie algebra $Lie(G)$: 
$$\eqalign{\langle a, b\rangle_{Lie(G)}&=  
tr_{Lie(G)}(ad_b^* ad_a )\cr &=  \sum_{ i=1}^d
\langle [a, c_i], [b, c_i]\rangle_{Lie(G)}\cr &
=-tr_{Lie(G)} [a, [b, \cdot]]\cr &= -tr_{Lie(G)}
[b, [a,
\cdot]].\cr} $$  where $c_i, i=1, \cdots, d$
varies in an O.N.B of $Lie(G)$.
 \mp   {\bf Proposition 7 }  \sp   {\it  
 The map: 
$$\eqalign{ ad: \Ci(S^1, Lie (G)) &\to
\Ci(S^1, Hom(Lie(G)))\cr
X:= \sum_{n\in \Z} a_n z^n &\to ad_X =\sum_{n\in
\Z}
ad_{ a_n} z^n\cr}$$ extends to a map: 
$$ad: H^{1\over 2}(S^1, Lie(G)) \to
g_{res}^D $$ with the notations of
(3.3) replacing $M$ by $S^1$ and $E$ by the
trivial bundle $S^1\times Lie(G)$. In particular,
 $ [ad_X]_{+-}$ and $[ad_X]_{-+}$    are
Hilbert-Schmidt operators on 
$ L^2(S^1, Lie(G)) $. }
\sp {\bf Proof: } (this proof is close to that of
proposition 6.3.1 in chapter 6 of [PS])
\sp For $b\in L$ , $q\in \Z$ and $X=
\sum_{n\in \Z} a_n z^n\in \Ci(S^1, Lie(G))$ we
have: 
$$\eqalign{ ad_X  (bz^q)&= \sum_{n\in \Z}
ad_{a_n} (b) z^{n+q}\cr
&= \sum_{p\in \Z} ad_{a_{p-q}} (b) z^p \cr}$$ so
that $ad_X$  is represented by an
$\Z\times \Z$ matrix with entries 
${ad_X}_{p, q}= ad_{a_{p-q}}\in Hom(Lie(G))$. 
  We compute the Hilbert-Schmidt norm of
${(ad_X}_{+-} $ denoted by $\Vert \cdot
\Vert_{HS}$. We shall use the fact that for $a
\in Lie(G)$ we have $\Vert ad_a \Vert_{HS}^2=
\Vert a\Vert^2$ where this last norm is the one
corresponding to the Killing form. 
$$\eqalign{ \Vert {ad_X}_{+-}\Vert_{HS}^2
&= \sum_{  q<0}\sum_{ i=1}^{dim Lie(G)} \langle
{(ad_X)}_{+-}b_iz^q,{(ad_X)}_{+-}b_iz^q\rangle 
\cr
&=  \sum_{p\geq 0, q<0} \sum_{ i=1}^{dim Lie(G)}
\langle
  ad_{a_{p-q}} 
(b_i)z^p, ad_{a_{p-q}} (b_i)z^p\rangle 
\cr
&=  \sum_{p\geq 0, q<0} \sum_{ i=1}^{dim Lie(G)}
\langle
  ad_{a_{p-q}} 
(b_i) , ad_{a_{p-q}} (b_i) \rangle_{Lie(G)} 
\cr
&= 
\sum_{p\geq 0, q<0} \Vert  {ad_X}_{p,
q} \Vert^2_{HS}  \cr
&= \sum_{k>0}  k \Vert   {ad_a}_{k}
\Vert^2_{HS}  \cr
&= \sum_{k>0} k\Vert    a_{k}
\Vert^2    \cr}$$ 
where $(b_i)_{i=1, \cdots, dim Lie(G))}$ is an
orthonormal basis of
$Lie(G)$. In the same way, we have: 
 $$  \Vert {ad_X}_{ -+}\Vert_{HS}^2 =  
\sum_{k<0} \Vert   a_{k} \Vert^2 
\vert  k\vert  .$$  
  From the above expressions  of the
Hilbert-Schmidt norms of ${ad_X}_{+-} $ and
${ad_X}_{ -+}  $, it follows  that these
are Hilbert-Schmidt whenever
$\sum_{k\in \Z} \vert k\vert  \Vert a_k
\Vert^2 $ is finite.  But this is   the
condition for $X$ to  lie in $H^{1\over 2}(S^1,
L)$. 
       \mp 
$\bullet$ {\it The first Chern form in terms of
the K\"ahler form}\sp
Recall that the exterior derivative of a form
can be expressed in terms of a Wodzicki residue
(see (5.8)).  Unlike in the finite dimensional
setting, here the
 weighted  first Chern form $r_1^Q$ might therefore
not be closed.  The following proposition shows
that it relates to two closed forms, the K\"ahler
form $\omega$ on the one hand (this fact had
already been shown by Freed) and the pull-back by
the adjoint representation of the cocycle $\l^D$ 
on the other hand.
\mp {\bf Lemma 9} \sp 
{\it The pull-back $ad^* \l^D$ of the cochain
$\l^D$ defined by (4.1) is closed on $H^{1\over
2}(S^1, Lie(G))$. } 
\sp {\bf Proof:} By Proposition 7, we
know that 
$ad$ takes its values in
$g^D_{res}$. On the other hand   $ \l_D$ is
indeed a cocycle on $g^D_{res}$. Combining these
two facts yields, for three left invariant vector
fields $X, Y, Z$: 
$$\eqalign{ \delta (ad^*\l^D)(X, Y, Z)&=ad^*\l^D( [X,
Y], Z)- ad^* \l^D ([Y, Z], X)+ ad^* \l^D([Z,
X] , Y)\cr &= \delta
\l^D ( ad_X, ad_Y, ad_Z)\cr &= 0\cr}$$
where $\delta\l^D$ denotes the coboundary of
$\l^D$. 
\mp 
\noi {\bf  Proposition 8}\
\sp
\noi  {\it  The first Chern form is given by:
$$r_1^Q(U, \bar V)= ad^*\l^{D}( U, \bar
V )=  -i
\omega(U,
\bar V) $$ 
for  $U, V\in H^{1\over 2}(S^1, Lie(G))$ and it is
a closed form. Here   $\omega$ is the symplectic
form on $H_e^{1\over 2} (S^1, G)$  defined in
section 8. }
\mp {\bf Proof:}  
From the results of Proposition 7 combined with
the description of $\l^D$ given in Lemma 3,  we
know that  
$\l^{D}(ad_U, ad_{\bar Y})=  tr( {ad_{\bar
V}}_{+-} {ad_U}_{-+} - {ad_U}_{+-} {ad_{\bar
V}}_{-+})$   is  an ordinary trace and
 hence it reads:
  $$ \l^{D}(ad_U, ad_{\bar V}) =  
 \sum_{q\in \N ,  i=1,\cdots, d} \left( \langle  
{ad_{\bar V}}_{+-} {ad_U}_{-+} z^q c_i, z^q
c_i\rangle  -   
\langle {ad_U}_{+-} {ad_{\bar V}}_{-+}   z^q c_i,
z^q c_i\rangle \right) 
 $$ 
where $c_i, i=1, \cdots, d$ varies in an O.N.B of 
$Lie(G)$. Since $L^2(S^1, Lie(G))$ is spanned  by
elements of the type $z^n a, n\in \Z, a\in
Lie(G)$, it is sufficient to compute this last
sum  for
$U=z^na$,
$\bar V=z^{-p} b$, $a, b\in Lie(G)$,
$n\in
\Z, p\in \Z$.  For $q\in \N, c\in Lie(G)$ we have: 
$$(ad_U)_{-+} (z^q c) = \matrix {  0 &\quad 
\hbox{ if } \quad n+q>0 \cr
                                     [a, c] z^{n+q} &\quad \hbox{ if }\quad n+q \leq 0\cr}$$ 
and
$$(ad_{\bar V})_{+-}(ad_U)_{-+} (z^q c) = \matrix
{  0 &\quad \hbox{ if }\quad n+q>0
\quad\hbox{or}\quad  n-p+q <0 \cr
[b, [a, c]]   
z^{n-p+q}  &\quad \hbox{ if }\quad n+q  \leq 0
\quad\hbox{ and}\quad  n-p+q \geq 0.\cr}$$ In the
same way 
$$(ad_{\bar V})_{-+} (z^q c) = \matrix {  0 
&\quad 
\hbox{ if }\quad -p+q>0 \cr
                                     [b, c] z^{n+q}  &\quad \hbox{ if } p+q \leq 0\cr}$$
and
$$(ad_U)_{+-}(ad_V)_{-+} (z^q c) = \matrix{  0 
&\quad  \hbox{ if }\quad -p+q>0\quad \hbox{or}
\quad n-p+q <0 \cr
                                     [a, [b, c]] 
z^{n-p+q}  &\quad  \hbox{ if }\quad -p+q  \leq 0
\quad\hbox{ and} \quad n-p+q \geq 0.\cr}$$ Finally
this yields (since the only non vanishing terms
correspond to $n-p=0$): 
$$ \eqalign{   \l^D(ad_U, ad_{\bar V}) &=
-\sum_{0< q\leq n,  i=1,\cdots, d} \langle [a, [b,
c_i]] z^q, c_iz^q\rangle\cr &= -n tr_{Lie
(G)}(c\to [a, [b, c]]) \cr &= n\langle
a,b\rangle_{Lie
(G)}\cr
 &=-i \omega(z^n a, z^{-n}b)\cr
 &=-i \omega(U, \bar V)\cr}$$  which yields the result. 
\mp A similar computation [F1]  (see Theorem
2.20)   shows that  the weighted first Chern form
can be expressed  in terms  of the symplectic form
$\omega$  and yields the
result. 
 \vfill
\eject
\noindent 
 \mp {\bf   Appendix A.  Classical   elliptic
pseudo-differential operator}\mp {\it This
Appendix 
 gives  a brief presentation of the basic
tools in our framework, namely classical
pseudo-differential operators and particularly
elliptic ones, their logarithms and their complex
powers. Classical references are [G], [Se],
[Sh].}  \mp  $\bullet${\it The symbol set}\sp  
Let $U$ be an open subset of
$\R^d$. Given $\alpha\in
\C$, let us denote by $S^\alpha(U   )$ the set
of  complex   valued  smooth
function
$$\eqalign{U \times \R^d &\to \R \cr
(U, \xi)&\to \sigma (U,
\xi)\cr}$$  satisfiying  the following property.
Given any  compact subset
$K$ of
$ U $ and any two multiindices $\gamma=
(\gamma_1, \cdots, \gamma_d)\in \N^d$, $\delta=
(\delta_1, \cdots, \delta_d)\in \N^d$ , there exists a
constant
$ C_{\alpha,
\beta}(K)$ such that
$$\vert D_x^\gamma D_\xi^\delta\sigma(x, \xi)\vert
\leq  C_{\gamma, \delta} (K) (1+ \vert
\xi\vert)^{Re\alpha-\vert \delta\vert}\quad \forall x\in K, \forall \xi \in \R^d.$$
An element of
$S^\alpha(x )$  is called a {\it symbol of order $\alpha$}. Let $S^{\leq m}(x)$
denote the set of symbols of order $\leq m$.
 The {\it principal part} of the symbol $\sigma
\in S^\alpha(x )$  ({\it or principal symbol }) is defined as follows:
$$\sigma_\alpha(x, \xi)= \lim_{t \to +\infty}
{\sigma(U, t\xi)\over t^\alpha}.$$
 \sp\noi A {\it smoothing symbol} is a symbol in
$$S^{-\infty } (x )\equiv
\bigcap_{k\in
\N} S^{\leq -k}(U) $$  and the relation
$$\sigma\simeq \tilde
\sigma\Leftrightarrow \sigma-\tilde \sigma \in
S^{-\infty } (U )$$
defines an equivalence relation on   $S (x )$. \sp
A symbol of order
$\alpha$ is called a {\it classical symbol} if there
exist $\sigma_{\alpha-j} \in  S^{\alpha-j}(x )$, $j\in \N$  such that:
 $$\sigma(x, \xi)\simeq \sum_{j=0}^\infty
\sigma_{\alpha-j}(x, \xi)  $$  which  are positively
homogeneous, i.e
  $$\sigma_{\alpha-j} (x, t\xi)= t^{\alpha-j}\sigma_{\alpha-j}(x,
\xi)\quad \forall t\in \R^{+ }.$$
\sp
Following Kontsevich and Vishik [KV 1], we shall
say that a classical symbol lies in the
  {\it odd-class } if the
positively homogeneous components $\sigma_{\alpha-j}$ are
moreover homogeneous i.e:
 $$\sigma_{\alpha-j}(x, t\xi)= t^{\alpha-j}\sigma_{\alpha-j}(x,
\xi)\quad \forall t\in \R.$$
\mp $\bullet$ {\it From symbols
to pseudo-differential operators}\sp
To a symbol $\sigma\in S^\alpha(U )$ we can
associate the {\it pseudo-differential operator (P.D.O.)
with symbol $\sigma$} defined by:
$$\eqalign{ A: C_c^\infty(x)&\to \C^\infty(\R^d)
\cr
u&\to  \left(U\to Au (x)  =
\int_{\R^d} e^{i\xi \cdot  xU }\sigma( U ,
\xi) \hat   u (\xi)  d\xi\right).\cr} $$
where $C_c^\infty(U)$ denotes the space of
complex
  valued  smooth functions with compact
support in
$U$. The {\it principal symbol} of $A$ is given by the principal part
$\sigma_P(A)$  of its
symbol
 $\sigma(A)$. If the symbol is classical,
we shall call the corresponding P.D.O
{\it classical}.
 The set of classical P.D.Os of order $\alpha$ (resp. $\leq m$) is denoted by
$PDO^\alpha(U)$ (resp. $PDO^{\leq m} (U)$).
\sp A {\it smoothing P.D.O. } is an operator   in
$$PDO^{-\infty } (U )\equiv
\bigcap_{k\in
\N} PDO^{\leq -k}(U ) $$ and there is an exact sequence:
$$0\to PDO^{-\infty } (U )\to PDO^{\leq m}(U ) \to S^{\leq m}(U )\to 0. $$
\sp An {\it ordinary differential operator} of order
$m\in \N$ is defined by a polynomial symbol (the
polynomial being of order $m$) in
$\xi$:
 $$\sigma(x, \xi)= \sum_{j=0}^m a_k(x) \xi^k. $$
A differential operator is {\it local} i.e $u=0
\Rightarrow Au=0$. But a pseudo-differential
operator, because of the smearing produced by the
Fourier transform, is not local. It is only
pseudo-local i.e if $u$ is smooth on an open
set $U$ then
$Pu$ is also smooth on any open subset $V\subset
U$.
\sp The various classes of symbols introduced
previously induce corresponding classes of
pseudo-differential operators.  A  {\it
classical pseudo-differential} is a P.D.O such
that its symbol has components given by classical
symbols and   an  {\it   odd-class classical  P.D.O }
is a classical P.D.O such that its symbol has
components given by  symbols in the odd class.  Notice
that ordinary differential operators with integer order provide
examples of P.D.Os in the odd class.
 \mp  $\bullet$ {\it pseudo-differential operators
acting on sections of vector bundles}\sp The notion of pseudo-differential
operator
 can
be carried out to operators acting on sections of
vector bundles.  Let
$\pi_E:E\to M$,
$\pi_F:F\to M$ be two smooth
 vector bundles with rank $r_E$ and $r_F$
respectively based on a smooth
   manifold
$M$ with dimension $d$. An     operator $$P:
\Gamma(M, E)\to \Gamma(M,F)$$ acting from the space
$\Gamma(M, E)$ of smooth sections of
$E$ to the space $\Gamma(M, F)$ of smooth sections of
$F$ is called a {\it pseudo-differential operator
of order $\alpha$ } if  given   a neighborhood   of any point $m\in M$,
there is a local trivialization i.e
a morphism:
$$\phi: (U, U\times \C^{r_E}, U\times \C^{r_F})\to (M, E,F)$$
where $U$ is an open subset of $\R^d$,
 the induced
linear map $ \phi^\sharp A: C_c^\infty  (U, U\times \C^{r_E})\to  C_c^\infty
(U, U\times \C^{r_F})$
has  a symbol
 $\sigma(\phi^\sharp A) \in C^\infty( U\times \R^d) \otimes {\cal M}_{r_E,
r_F} (\C)$
   with matrix components in $S^\alpha(U)$. Here
${\cal M}_{k, l}$ denotes the space of $k\times l$ matrices with
coefficients in $\C$.
It is a {\it classical P.D.O } if  $\sigma(\phi^\sharp A) $ is a classical
symbol.
These definitions involve a choice of trivialization but can be shown to
 be independent of this choice.
\sp $\sigma(\phi^\sharp A) $ is called
 the {\it (formal) symbol} of $A$ and  is only
defined locally. However   its principal part $\sigma_P(\phi^\sharp A) $
is independent of the choice of coordinate charts
and is therefore defined globally. We shall denote it by $\sigma_P(A).$ It
is  called the
{\it principal symbol } of the P.D.O.
\sp  Let us denote by $ PDO^\alpha (M, E,F)$   the space of all classical
 P.D.Os of order
  $\alpha$ and by    $ PDO^{\leq m} (M, E,F)$ the space of all classical
 P.D.Os  of order
  $\leq m $.  When $E=F$, we shall denote these spaces by $PDO^\alpha (M, E)$
(resp.
$PDO^{\leq m}  (M, E)$ )
and when $E$ is the trivial bundle $M\times \C$
 by $PDO^\alpha(M)$ (resp.$PDO^{\leq m}  (M )$).
The {\it symbol set} $S^{\leq m} (M, E,F)$ is defined by the exact sequence:
 $$0\to PDO^{-\infty } (M, E, F) \to PDO^{\leq m} (M, E, F) \to
S^{\leq m}(M, E, F) \to 0  $$
where $PDO^{-\infty } (M, E, F) := \bigcap_{k\geq 0} PDO^{\leq -k} (M, E,F)$.
  \sp When $M$ is compact, their is a notion of product of
two pseudo-differential operators and
 $$PDO(M, E):= \bigcup_{m\in \Z} PDO^{\leq m} (M, E)$$
defines an associative algebra.
 \mp
{\it From now on we shall assume $M$ is a smooth compact manifold without
boundary}.
\mp
$\bullet${\it  Admissible elliptic pseudo-differential
operators}\mp
Let $E$ and $F$ be as before two finite rank
vector bundles based on a smooth manifold $M$.
When
$r_E=r_F$  a pseudo-differential operator
$P:
\Gamma(E)\to
\Gamma(F)$ of order $m$ is called {\it elliptic  }
if its principal symbol $\sigma_P^m(x, \xi)$ is
an invertible matrix for $\xi\neq 0$. Let
$Ell(M, E)$ denote the set of elliptic classical
P.D.Os.
 \sp   Let
us denote by
$Ell_{ord>0}^*(M, E)$ the class of
 invertible elliptic operators of strictly
positive order. Since $M$ is compact the spectrum
$spec(A)$ of such an  operator consists of
isolated eigenvalues with finite multiplicity [Sh].
There is therefore a disc $D_R$ of radius $R>0$
around the  origin which does not contain any
point of the spectrum.
 We shall say that $A$ has {\it   spectral cut}
$L_\theta$ if there is   a   ray
$ L_\theta= \{\l\in \C, arg \l=\theta\}$ in the complex plane which
does not intersect the spectrum of $A$.  Such an operator will be called
{\it admissible} and we shall denote by
$Ell_{ord>0}^{*,adm}(M, E)$ the set of such admissible operators.
Any invertible elliptic operator with  strictly positive order such that
the  matrix given by its principal symbol has no eigenvalues
in some non empty conical
 neighborhood $\Lambda$ of a ray in the spectral plane is
admissible since  in that case at most a finite number
of eigenvalues of the operator are
 contained in  $\Lambda$ [Sh]. \sp 
  Let us introduce some notations. 
$Ell^{ s.a}_{ord >0}(M, E)$, resp. 
 $Ell^{ +}_{ord >0}(M, E)$ denotes 
the set of self-adjoint, resp. positive   self-adjoint elliptic operators with
 strictly positive order.
 Adding an upper index $*$ restricts to injective operators so that
 $Ell^{ *s.a}_{ord >0}(M, E)$, resp. 
  $Ell^{*+}_{ord >0}(M, E)$ denotes 
the set of self-adjoint injective, resp. positive self-adjoint injective  elliptic operators with
 strictly positive order and  we have following inclusions:
  $$Ell^{*+ }_{ord>0}(M, E)\subset Ell^{* s.a}_{ord >0}(M, E) \subset 
 Ell_{ord>0}^{*,adm}(M, E).$$
     \mp $\bullet$ {\it Complex powers and logarithms of elliptic operators}\sp
Let $A\in Ell^{*,adm}_{>0}(M, E)$ with spectral cut $L_\theta$.
 For
$Rez < 0$, the complex power
$A_\theta^z$ of $A$ is a bounded operator on any
space
$H^s(M, E)$ of  sections of $E$ of Sobolev class
$H^s$  defined by the contour integral:
 $$A_\theta^z={i\over 2\pi} \int_{\Gamma_\theta}
\l^z (A-\l I)^{-1} d\l$$
where $\Gamma_\theta= \Gamma_{1, \theta}\cup
\Gamma_{2,
\theta}\cup \Gamma_{3, \theta}$
$\Gamma_{1,\theta}= \{\l=  r e^{i\theta}, r\geq
R\}$,$\Gamma_{2,\theta}= \{\l= R
e^{i \phi }, \theta\geq\phi\geq -\theta
 \}$, $\Gamma_{3,\theta}= \{\l=  r
e^{ i (\theta-2\pi)  }, r\geq R\}$.
 Here $\l^z= exp (zlog \l)$
where $log \l= log \vert \l\vert + i \theta$ on
$\Gamma_{1, \theta}$ and $log \l= log \vert
\l\vert  + i (\theta-2\pi)$ on
$\Gamma_{3, \theta}$. \sp
This definition is independent of the choice of
 $R$ but depends on the choice of $\theta$ and
yields for any
$z\in
\C$  an elliptic operator $A_\theta^z$ of order
$z\cdot ord(A)$. When
$z=-k, k\in
\N$, then
$A^z$ coincides with $A^{-k}$ of order
$-k\cdot ord(A)$. When $M$ is Riemannian, $E$
is hermitian and $A$ is essentially self-adjoint,
then $A_\theta^z$ is independent of the choice of
$\theta$ and coincides with the complex powers
defined using spectral representation.
\sp {\it  In the following we shall focus on operators in
$Ell^{ +}_{ord>0}(M, E)$ in which case
we shall use the principal branch  of the logarithm, taking
$\theta=\pi$ and simply drop the mention $\theta$.}
\mp
For arbitrary $k\in \Z$, the map  $z\to
A^z_\theta$  defines a holomorphic function from
$\{z\in \C, Rez<k\}$ to the space of bounded linear maps ${\cal L} (H^s(M, E)\to
H^{s-kord A}(M, E))$ for any $s\in \R$ and we can
set:
$$log_\theta A\equiv \left[{\partial \over
\partial z} A_\theta^z\right]_{z=0} $$
which defines a (non classical) P.D.O operator of zero order and hence a
bounded operator from
$H^s(M, E) $ to $H^{s-\e}(M, E)$ for any $\e>0$
and any $s\in \R$.  In local coordinates $(x,
\xi) $ on $T^*M$, the  symbol of $log_\theta A$
reads:
$$\sigma_{log_\theta A}(x, \xi)= ord(A) log \vert
\xi\vert Id +
 \hbox{ a classical P.D.O symbol of order 0}.$$
Hence, although the logarithm of an injective
elliptic classical pseudo-differential operator
with admissible cut
$L_\theta$ is not itself a classical
pseudo-differential operator, for two operators
$A\in Ell_{ord>0}^{*, adm} (M, E)$, $B\in Ell_{ord>0}^{*, adm}(M,
E)$ admitting spectral cuts $L_\theta$ and
$L_\phi$:
$${log_\theta A\over ord A }-{log_\phi B\over ord B }\in
PDO^0(M, E). 
$$ \vfill \eject \centerline{\bf Appendix B}

\mp{\it  In this appendix, we recall why the first fundamental group 
of $Ell_{0 }^{*,0}(M,E)$ is non trivial, where $Ell_{0 }^{*,0}(M, E) $ denotes the pathwise connected
component of identity in the group of invertible zero order elliptic P.D.Os.}
 \bp  $E$ denotes a finite rank vector bundle based on a compact boundaryless Riemannian manifold $M$
equipped with a connection $\nabla$.
We keep the notations of Appendix A.    To begin with, let us describe the  topology on
$Ell_0^*(M,E)$.  
\mp   $\bullet${\it A Fr\'echet structure on  the algebra of P.D.Os and their symbols} 
 \sp 
The space $SPDO(M, E)$ of symbols of classical pseudo-differential operators is a Fr\'echet space
when equipped with the following family of semi-norms  labelled
 by multiindices
$\gamma
\in
\N^d$, $\delta \in \N^d$ and $i\in \{1, \cdots, N\}$, $k\in \N$: 
$$ \|\sigma \|_{\gamma,\delta,k}:=max_{i} \Big( sup_{y \in \overline{ V_i }, 
\xi \in \R^d, \Vert \xi \Vert=1 }   \| D_y^\gamma D_\xi^\delta\sigma_k (y, \xi) \|
\Big)
$$ where     $\sigma_k$ is the homogeneous component of order $k$, $ \{V_i\}_{i=1, \cdots, N } $  is
a  finite open cover (with $\bar V_i$ compact) of
$M$ associated to a partition of unity $ \{V_i, \xi_i\}_{i=1, \cdots, N }  $ on $M$ subordinated to
some finite open covering 
$
\{U_i\}_{i=1, \cdots, N }  $.
\sp This topology combined with natural semi-norms on kernels of compact operators,   induces a
Fr\'echet structure  on the  space of  classical pseudo-differential
 operators  
$PDO (M,E)$ via the identification of operators with their symbol up to a smoothing operator (see  
[KV1] section 3)
 \mp    
 $\bullet${\it  From elliptic operators to their principal symbols} 
\sp Let $S_P  Ell_0^*(M,E) $  denote the group of principal symbols of operators in
$Ell_0^*(M,E)$ and let  $ Ell_0^{*,0}(M,E) $, resp.  $S_PEll_{0 }^{*,0}(M, E)$ be the pathwise
connected component of identity of these topological spaces.   
  \mp{\bf Lemma} \sp                                    
  $$\Pi_1(S_P  Ell_{0 }^{*,0}(M,E) )= \Pi_1(Ell_{0 }^{*,0}(M,E)).$$
 \sp {\bf Proof:} We show that $Ell_{0 }^{*,0}(M,E)$ and $S_P  Ell_{0 }^{*,0}(M,E) $ have same
homotopy type. The result then follows. 
\sp  Let $$\eqalign{ f:  Ell_{0 }^{*,0}(M,E) &\to S_P 
Ell_{0 }^{*,0}(M,E) \cr
 A&\to \sigma_P(A).\cr}$$ 
The connection on  $E$ and the metric on $M$ give   a canonical way of assigning to a
symbol $\sigma$, an operator $Op(\sigma (A))(u)(x)$ with same
principal symbol (see e.g. [BML] page 188 and references therein):
$$Op(\sigma (A))(u)(x):=  \int_{T_x^*M}e^{i  x \cdot \xi} \sigma (\xi) \hat u(\xi)d\xi  $$ 
  where the Fourier transform is defined using the exponential map on $M$ and parallel transport on
$E$. This gives rise to a map:
$$\eqalign{ g: S_P  Ell_{0 }^{*,0}(M,E)  &\to  
Ell_{0 }^{*,0}(M,E)\cr
 \sigma&\to Op\left(\sigma \right).\cr}$$ 
 \sp
Then $f\circ g= Id$. One does not expect that $g\circ f= Id$ because  $A = Op(\sigma_P(A))( 1+ K)$
for some (uniquely defined) compact operator $K$.  However we have $g\circ f\sim  Id$ where $\sim$
denotes homotopy of maps. Indeed, we can build a map 
$(A, s)\to F(A, s) $ on $Ell_0^*(M,E)\times [0, 1]$ such that $F(A, 0)= g \circ f(A)=
Op\left(\sigma_P(A)\right)$ and $F(A, 1)= A$ for any $A \in Ell_0^*(M,E)$. Set $$F(A,
s):=(1-s )Op\left(\sigma_P(A)\right)+ s A=Op\left(\sigma_P(A)\right) (1+sK) ,$$  then $F(A,
s)\in  Ell_0  (M,E)$ since $\sigma_P(F(A, s))=
\sigma_P\left(Op\left(\sigma_P(A)\right)\right)= \sigma_P(A)$  and
$A
$ is elliptic of order zero. The operator $1+s K$ might not be invertible so  let us show how we can
modify it continuously into  a family  
$\tilde F(A, s)\in  Ell_0^* (M,E)
$.  Let $\gamma(s):=Id+K_s$.  Let us 
 denote by $S_K $ the spectrum of $K $. Since $K $ can be extended into a compact operator of
$L^2(M,E)$, $S_K $ is compact and consists of isolated points, up to 0 which may be an
accumulation point. The operator $Id+sK $ is not invertible whenever 
  $-{1\over s}$ is an
eigenvalue of $K $. Since $s\in[0,1]$, we have $-{1\over s} \leq -1 $. The set $ ]-\infty,-1]
\cap S_K  $ is finite, so $ \gamma(s) $ is invertible, up to a finite number of $s\in [0,1]$.  
We can therefore modify  $\gamma$ continously into some continuous path $\tilde \gamma(s):= 
1+\phi(s)
K$    with values in the set
$
\{ Id+K
\ invertible : K 
\in PDO^{\leq-1}(M,E)
\}
$, where $PDO^{\leq \alpha}(M,E)$ denotes the subspace of pseudodifferential opertors with order
$\leq \alpha$.   Setting
$\tilde F(A, s):=  Op\left(\sigma_P(A)\right)\left(1+\phi(s)K\right)$ yields the result. 
 \mp $\bullet$ {\it The fundamental group     $
\Pi_1( Ell_0^*(M,E)) $}
 \sp We shall discuss the case when $E$ is trivial, but the results extend to the non trivial
case using the fact that $M$ being compact, for large enough $N$, 
there is a vector bundle $F$ based
on
$M$ such that $E\oplus F\simeq 1_N$ where $1_N$ is the trivial bundle of rank $N$ on $M$.
\sp Let
${\cal A}$ be a complex unitary Banach algebra and let
$N$ be a positive integer large enough. Then ([Ka], chIII, Th1.25) the fundamental group
$\Pi_1(GL_N({\cal A}))$ is isomorphic to
$K({\cal A})$, where 
$GL_N({\cal A})$ is the group of invertible maps of ${\cal A}^N$ and $K({\cal A})$ is
 the Grothendieck group of finitely generated projective ${\cal A}$-modules. Applying that result
to the Banach space $C(S^*M)$ of continuous functions on $S^*M$, we find that:
$$ \Pi_1(GL_N(C(S^*M))) \simeq  K_0(S^*M):= K(C(S^*M)).$$
 \sp Using the classical result (see e.g. [Di], vol. 3) that on a compact manifold a continuous
function is homotopic to a smooth function and that two homotopic smooth functions are smoothly
homotopic, we have:
$$ \Pi_1(GL_N(C^{\infty}(S^*M)))=\Pi_1(GL_N(C(S^*M))). $$
It follows   that for the trivial
bundle $1_N$  based on $M$ of  rank $N$ large enough: 
 $$\Pi_1\left(S_P
Ell_0^{*,0}(M,1_N)\right)\simeq \Pi_1(GL_N(C(S^*M)))\simeq  K_0(S^*M).$$
\sp On the other hand, $\Pi_1(GL_N(C(S^*M))) $  being generated by one parameter subgroups $exp(2\pi
i t p)$,
$p:\left( C(S^*M)\right)^N\to \left(C(S^*M)\right)^N  $, such that $p^2=p$ [Ka],  is non trivial so
that 
$
\Pi_1(  Ell_0^*(M,1_N))
$ is   non trivial as announced.
\sp {\it Remark:} In fact, since the arcwise connected component of identity $Ell^{ *,0}(M, E)$ of
the  group
$Ell^*(M, E)$ of invertible elliptic operators has the same fundamental group as  
$Ell^{*,0}_0(M, E)$  consisting of those that have zero order, this shows that 
$\Pi_1( Ell^{*,0}(M, E))= \Pi_1( Ell_0^{*,0}(M, E) )\simeq  K_0(S^*M) $ is generated by one parameter
subgroups $exp(2\pi i t P)$ of the above type and hence non trivial. 
  \mp \vfill
\eject
\noindent\centerline{\bf REFERENCES}
\vskip 0,5cm
\item{[AP]} M.Arnaudon, S.Paycha, {\it Regularisable and minimal orbits for
group actions in infinite dimensions}, Comm.Math.Phys. {\bf 191} , 641--662
(1998)
\item{[BF]} J.M.Bismut, D.Freed,{\it  The Analysis of elliptic families I},
Comm.Math.Phys. {\bf 106}, 159--176 (1986)
 \item{[BGV]} N. Berline, E. Getzler, M. Vergne, { \it Heat-kernels and
Dirac operators}, Springer Verlag 1991
\item{[CFNW]}M.Cederwall, G.Ferretti, B. Nilsson, A.Westerberg, {\it
Schwinger terms and Cohomology of pseudodifferential
operators}, Comm.Math.Phys. {\bf 175} , 203--220 (1996)
\item{[CE]} H. Cartan, S. Eilenberg, {\it Homological algebra}, Princeton
University Press (1956)\item{[C1]} A. Connes,{\it  Non commutative
differential geometry},   Publications  I.H.E.S , n.62 (1985) p. 257-360
\item{[C2]} A. Connes,{\it  Non commutative geometry},   Academic Press,
San Diego, 1994
\item{[D]} C.Ducourtioux,{\it Ph.D thesis in
preparation}
\item{[Di]} J. Dieudonn\'e, {\it El\'ements d'analyse 3}, Gauthier-Villars (1970)
\item{[DL]} B. Driver, T. Lohrenz, {\it Log
Sobolev Inequalities on loop groups }, Journal of Functional Analysis  {\bf 140}, 381-438 (1996)
\item{[FLPTT]} N. Fagella, A. Lesne, S. Paycha, L. Tedeschini-Lalli, S.T.
Tsou, {\it Renormalization, Proceedings of a workshop}, Paris 1997
\item{[F1]} D. Freed, {\it The geometry of loop groups}, Journ. Diff. Geom. {\bf 28},
  223--276 (1988)
\item{[F2]} D. Freed, {\it Flag manifolds and infinite dimensional K\"ahler geometry}, in "Infinite dimensional groups with applications, Ed. V. 
Kac, Math. Sci. Res. Inst. Publ., Vol 4, Springer, New York, 1985
\item{[G]} P. Gilkey, {\it Invariance
theory. The heat equation and the Atiyah-Singer index theorem},  Publish or
Perish (1984), second edition, CRC Press (1995) 
\item{[H]} L.H\"ormander, {\it The analysis
of linear partial differential operators}, I. Grundlagen der mathematischen
Wiss. {\bf 256}, Springer Verlag (1983)
\item{[Ka]} M.Karoubi, {\it $K$-theory (An introduction)}, Springer Verlag 1978 
 \item{[K]} Ch. Kassel, {\it Le r\'esidu non commutatif [d'apr\`es
Wodzicki]}, S\'eminaire Bourbaki {\bf 708} (1989)
\item{[Ki]}A. Kirillov, {\it Elements of theory of representation}, Berlin,
Springer Verlag, 1976\item{[KOS]}H.H Kuo, N. Obata, K. Saito, {\it L\'evy
Laplacians of
generalized functions on a nuclear space},
\item{[KK]} O.S. Kravenchko, B.A. Khesin, {\it A central extension of the
algebra of pseudodifferential symbols}, Funct.Anal.Appl. {\bf  25},
152--154  (1991)
\item{[KV1]} M.Kontsevich, S. Vishik, {\it Determinants of elliptic
pseudodifferential operators}, Max Planck Preprint (1994)
\item{[KV2]} M.Kontsevich, S. Vishik, {\it Geometry of determinants of
elliptic operators} in Functional Analysis on the Eve of the 21st Century
{\bf Vol. I }\hfill \break \noi (ed. S.Gindikin, J.Lepowski,R.L.Wilson)
Progress in Mathematics (1994)
  \item{[Ku]} N.H.Kuiper, {\it The homotopy type of
the unitary group of Hilbert spaces}, Topology 3
(1965)
 \item{[L]} S. Lang, {\it An introduction
to differentiable manifolds}, Interscience,  John Wiley (1967) 
\item{ [LM]} B.Lawson, M.L Michelsohn, {\it Spin Geometry, Princeton Mathematical Series}, 1989
\item{[Le]} M.Lesch, {\it On the non commutative
residue for pseudo-differential operators with
log-polyhomogeneous symbols}, Annals of
global analysis and geometry {\bf 17} (1998) 
\item{[MRT]} Y. Maeda,  S.
Rosenberg, P.Tondeur{\it The curvature of gauge orbits } in
" Global Analysis and Modern Mathematics" edited
by K. Uhlenbeck, Publich or Perish (1994)
 \item{[M]} J. Marsden, {\it Applications of Global Analysis in Mathematical Physics},
 Publish or
 Perish, 1974
 \item{[M1]} J.Mickelsson, {\it Current algebras and groups},New York,
Plenum Press, 1989
\item{[M2]} J.Mickelsson,
{\it Second Quantization, anomalies and group extensions}, Lecture notes
given at the
 "Colloque sur les M\'ethodes G\'eom\'etriques en physique,
C.I.R.M, Luminy, June 1997
\item{[MN]} R. Melrose and V. Nestor,
{\it Homology of pseudo-differential
operators I. manifolds with boundary},
Preprint: funct-an/9606005, Oct. 98\item{[Pa]} S.
Paycha, {\it Renormalized traces as a looking
glass into infinite dimensional geometry; a
proposal to extend concepts of Riemannian geometry
to a class of Hilbert manifolds}, Preprint 1999
 \item{[Pr]} A. Pressley, {\it The
Energy flow on, the loop space of a
compact Lie group}, Journ. of the
LOndon Math. Soc. Vol 26, 1982
\item{[PS]} A.
Pressley, G. Segal, {\it Loop groups}, Oxford
Univ.Press (1988)
\item{[R]} A.O.Radul, {\it Lie algebras of differential operators, their
central extensions, and W-algebras}, Funct.Anal.Appl. {\bf 25},
25--39(1991)
\item{[Rog]} C. Roger, {\it extensions centrales d'algebres et de groupes
de Lie de dimension
infinie, algebres de Virasoro et g\'en\'eralisations},
Pr\'epublications de l'institut Girard
 Desargues (1995)
\item{[S]} J. Schwinger, {\it Field theory commutators}, Phys. Rev. Lett.
{\bf 3}, 296-297 (1959)
\item{[Se]} R.T. Seeley, {\it Complex powers of an elliptic operator,
 Singular integrals}, Proc. Symp. Pure Math., Chicago, Amer.Math. Soc.,
 Providence,
 288-307   (1966)  
\item{[Sh]} M.A. Shubin, { \it pseudodifferential Operators and Spectral
Theory}, Springer Verlag 1987
 \item{[SW]} M. Spera, T. Wurzbacher, {\it 
Differential geometry of Grassmannian embedding of
based loop groups}, Preprint  1998
 \item{[W]} M. Wodzicki, {\it Non commutative residue} in  Lecture Notes in
Mathematics {\bf 1289} Springer Verlag (1987)
\item {[Wu]} T. Wurzbacher, {\it Symplectic
geometry of the loop space of a Riemannian
manifold}, Journal of Geometry and Physics {\bf
16}, 345-384 (1995) \end
\end \bye